\documentclass[a4paper,11pt,reqno]{amsart}

\usepackage{amssymb,enumitem,color,fancyvrb,prooftree}
\usepackage[colorlinks=true,urlcolor=blue,pagebackref=true]{hyperref}
\usepackage{threeparttable}
\usepackage[a4paper,margin=3cm,headsep=1cm,top=3cm]{geometry}
\usepackage{setspace}
\allowdisplaybreaks
\newtheorem{theorem}{Theorem}[section]

\newtheorem{lemma}[theorem]{Lemma}
\newtheorem{corollary}[theorem]{Corollary}

\theoremstyle{definition}
\newtheorem{definition}[theorem]{Definition}
\newtheorem{notation}[theorem]{Notation}

\newtheorem{remark}[theorem]{Remark}

\newcommand{\FO}{\mathbb{P}}
\newcommand{\Prod}{{\Pi}}
\newcommand{\Sum}{{\Sigma}}
\newcommand{\I}{\mathrm{I}}
\newcommand{\W}{\mathrm{W}}

\newcommand{\pair}[1]{\langle{#1}\rangle}
\newcommand{\V}{\mathrm{V}}
\newcommand{\va}{\mathrm{V_{ext}}(A)}
\newcommand{\vo}{\mathrm{V_{ext}}(\omega)}
\newcommand{\vset}[1]{\{{#1}\}_{A}}
\newcommand{\vpair}[1]{\langle{#1}\rangle_{\!A}}
\DeclareMathOperator{\dom}{\mathrm{dom}}
\DeclareMathOperator{\op}{\mathrm{OP}}     
\DeclareMathOperator{\ran}{\mathrm{ran}}

\DeclareMathOperator{\succe}{\mb{succ}}
\DeclareMathOperator{\pred}{\mb{pred}}
\DeclareMathOperator{\fun}{\mathrm{Fun}}

\newcommand{\fo}{\Vdash}   
\newcommand{\fof}{\Vdash_\FO}

\newcommand\mb{\mathbf}      
\newcommand{\vp}{\varphi}
\newcommand\dash{\text{-}}
\newcommand{\czf}{{\sf CZF}}
\newcommand{\izf}{{\sf IZF}}
\newcommand{\ac}{{\sf AC}}
\newcommand{\cac}{{\sf CAC}}
\newcommand{\dc}{{\sf DC}}
\newcommand\ft{\mathsf{FT}}
\newcommand{\rdc}{{\sf RDC}}
\newcommand{\pax}{{\sf PAx}}
\newcommand{\rea}{{\sf REA}}
\newcommand{\imp}{\rightarrow}

\newcommand{\biimp}{\leftrightarrow}

\def\rul#1#2#3{\prooftree     #1 \justifies #2 \using{#3} \endprooftree}
\title[Realizability for choice in dependent types]{Extensional realizability and choice for dependent types in intuitionistic set theory}
\thanks{The author's research was funded by the Alexander von Humboldt foundation.}

\author[Frittaion]{Emanuele Frittaion}
\address{Department of Mathematics, Technische Universit\"{a}t Darmstadt, Germany}


\begin{document}
	
	
\begin{abstract}	
In \cite{FR21}, we introduced an extensional variant of generic realizability \cite{M84}, where realizers act extensionally on realizers, and showed that this form of realizability provides {\em inner} models of $\czf$ (constructive Zermelo-Fraenkel set theory) and $\izf$ (intuitionistic Zermelo-Fraenkel set theory), that further validate $\ac_\ft$ (the axiom of choice in all finite types). In this paper, we show that  extensional generic realizability validates several choice principles for  {\em dependent types}, all exceeding $\ac_\ft$.   We then show that adding such choice principles does not change the arithmetic part of either $\czf$ or $\izf$. 
\end{abstract}
	
\maketitle

\section{Introduction}
In this paper, we prove that extensional generic realizability \cite{FR21} validates the natural extension of $\ac_\ft$ (the axiom of choice in finite types) and $\rdc_{\ft}$ (relativized dependent choice in finite types)    to  {\em transfinite} dependent types.
By  routinely  combining extensional generic realizability  with a  Beeson-style forcing construction, along the lines of \cite[Section 5]{friedman_scedrov84}, we also show that  $T$ augmented with such choice principles is conservative over $T$ for arithmetic sentences, for a host of intuitionistic set theories $T$. In particular, we extend  Friedman and Scedrov's result \cite{friedman_scedrov84} that  $\izf$ augmented with $\cac_{\ft}$ (finite type countable choice) and $\dc_{\ft}$ (finite type dependent choice) is arithmetically conservative over $\izf$.\footnote{$\cac_{\ft}$ is the schema $\forall n\, \exists y^{\tau}\, \vp(n,y)\imp \exists f^{0\to\tau}\, \forall n\,\vp(n,f(n))$, while $\dc_{\ft}$ stands for $\forall x^{\sigma}\,\exists y^{\sigma}\,\vp(x,y)\imp \forall x^{\sigma}\, \exists f^{0\to\sigma}\, (f(0)=x\land\forall n\,\vp(f(n),f(n+1)))$.}
	
Generic realizability\footnote{The name ``generic'' for this kind of realizability is due to  McCarty \cite[p.\ 31]{M84}.} is distinguished by its treatment of unbounded quantifiers: a realizer of $\forall x\, \vp(x)$ (resp.\ $\exists x\, \vp(x)$) must be a realizer of $\vp(x)$ for every   $x$ (resp.\ some  $x$).  
This strain of realizability  goes back to Kreisel and Troelstra's  realizability for theories of elementary analysis with species variables  \cite[Section 3.7]{KTr70}. A direct descendant of Kreisel and Troelstra's  realizability  was  applied  to  nonextensional theories of higher order arithmetic  and  set theory by Friedman \cite{F73} and  Beeson \cite{B79,B85}.  McCarty \cite{M84,M86} developed a version of generic realizability that applies directly to the extensional set theory  $\izf$. He  also employed,  inspired by  Feferman  \cite{F75},   arbitrary  partial combinatory algebras, rather than just natural numbers. Whereas McCarty's approach  is  geared towards $\izf$, it was shown in \cite{R06}  that $\czf$ suffices for a formalization of generic realizability, thus providing a self-validating semantics for $\czf$ and other fragments of $\izf$.  
	
In \cite{FR21}, we added an extra layer of extensionality to McCarty's  generic realizability, and thus obtained    {\em inner} models of both $\czf$ and $\izf$ that further validate   $\ac_{\ft}$.\footnote{Extensional notions of realizability have been mainly considered in the context of arithmetic and finite type arithmetic (see Troelstra  \cite{T98} and van Oosten  \cite{ Oosten97,Oosten08} for more information, and \cite{F19, BS18} for recent applications). Versions of extensional realizability in the context of both arithmetic and set theory can be found in Gordeev  \cite{gordeev}.} 
In the present paper, we  show that  extensional generic realizability yields an interpretation of 
\begin{itemize}[leftmargin=5mm]
		\item  $T+\Prod\Sum\I\dash\ac^\dagger+\Prod\Sum\I\dash\rdc^\dagger$ in $T$, and
		\item  $T+\rea+ \Prod\Sum\W\I\dash\ac^\dagger +\Prod\Sum\W\I\dash\rdc^\dagger$ in $T+\rea$,
\end{itemize}
where $T$ can be either  $\czf$, $\izf$ or an extension thereof with other choice axioms. We will introduce this family of choice principles for dependent types in Section \ref{sec:prel}. Suffice it  to say here that  $\Prod\Sum\I\dash\ac^\dagger$ and $\Prod\Sum\W\I\dash\ac^\dagger$  are closely related to, and in fact follow from,   Aczel's $\Prod\Sum\I\dash\ac$ and $\Prod\Sum\W\I\dash\ac$. The latter are validated by the type-theoretic  interpretations of $\czf$ \cite{A78, A82} and $\czf+\rea$ \cite{A86}, respectively. Our interpretation does not appear to validate   $\Prod\Sum\I\dash\ac$ and relatives, even  if we assume it in the background theory.  We will see in Section \ref{sec: remarks} that extensional generic realizability refutes the statement that  every set is an image of a $\Prod\Sum\I$-set ($\Prod\Sum\W\I$-set), which instead holds true in Aczel's type-theoretic interpretation of $\czf$ ($\czf+\rea$).
	
The basic plan for  interpreting choice via extensional generic realizability can be easily described as follows. Consider the $\Prod\Sum\I$ case. We start out with a realizability universe $\va$, which is built on top of a given partial combinatory algebra $A$. We then construct  an {\em extensional} type structure on $A$ (Definition \ref{types}).  In particular, types are elements of $A$.  Next, we single out  {\em nice} sets of $\va$ of the form $X_\sigma$, where $\sigma$ runs over all types (Definition \ref{nice}). We think of these sets as the extensions of the corresponding types. With this identification in mind, we  show that $\ac_{\sigma,\tau}$ and $\rdc_\sigma$ hold in $\va$ (Theorem \ref{canonical choice} and Theorem \ref{canonical dependent}). Finally, we prove that every $\Prod\Sum\I$-set  of $\va$ is {\em nicely} represented (Theorem \ref{realizers} and Theorem \ref{canonical}).   Niceness is captured by a suitable notion of being injectively presented in $\va$ (similar notions can be encountered in \cite{A82,A86,R05a,R06b,R03m}). It will be crucial to work with extensional type structures.

\section{Preliminaries}\label{sec:prel}
	
\subsection{$\czf$ and $\izf$}
In this paper, we will be concerned with  (extensions of)  $\czf$ (constructive Zermelo-Fraenkel set theory) and $\izf$ (intuitionistic Zermelo-Fraenkel set theory). The logic is  intuitionistic first order logic with equality in the language of classical $\sf ZF$. As usual, a formula is bounded, or $\Delta_0$, if all quantifiers appear in the form $\forall x\in y$ and $\exists x\in y$.  \vspace{1mm}
	
The theory $\czf$ consists of the following axioms: \vspace{1mm}
	
\begin{threeparttable}
\begin{tabular}{ll}
			\textbf{Extensionality} &  $\forall z\, (z\in x\biimp z\in y)\imp x=y$ \\
			\textbf{Pairing} & $\exists z\, (x\in z\land y\in z)$\\
			\textbf{Union}  &  $\exists y\, \forall u\, \forall z\, (u\in z\land z\in x\imp u\in y)$\\
			\textbf{Infinity} &  $\exists x\, \forall y\, (y\in x\biimp y=0\lor \exists z\in x\, (y=z\cup\{z\}))$\tnote{$\ast$}\\ 
			\textbf{Induction} &  $\forall x\, (\forall y\in x\, \vp(y)\imp \vp(x))\imp \forall x\, \vp(x)$ 
			\	for all formulae $\vp$  \\
			\textbf{Bounded Separation} &  $ \exists y\, \forall z\, (z\in y\biimp z\in x\land \vp(z))$ 
			\  for  $\vp$ bounded \\
			\textbf{Strong Collection} &    $\forall u\in x\, \exists v\, \vp(u,v)\imp \exists y\, \psi(x,y)$ 
			\ for all formulae  $\vp$ \\
			& where $\psi(x,y)$ stands for \\
			&  $\forall u\in x\, \exists v\in y\, \vp(u,v)\land \forall v\in y\, \exists u\in x\, \vp(u,v)$\\
			\textbf{Subset Collection} & $\exists z\, \forall r\, (\forall u\in x\, \exists v\in y\, \vp(u,v,r)\imp \exists  y_0\in z\, \psi(x,y_0,r))$ \\
			&  for all  formulae $\vp$ 
			where  $\psi(x,y_0,r)$ stands for \\
			&  $\forall u\in x\, \exists v\in y_0\, \vp(u,v,r)\land \forall v\in y_0\, \exists u\in x\, \vp(u,v,r)$	
\end{tabular}
\begin{tablenotes}\footnotesize
			\item[$\ast$] Let $x=0$ be $\forall y\in x\, (y\notin x)$ and $x=y\cup\{y\}$ be $\forall z\in x\, (z\in y\lor z=y)\land \forall z\in y\, (z\in x)\land y\in x$.
\end{tablenotes}
\end{threeparttable}	
	
\medskip

The theory $\izf$  consists of  extensionality, pairing, union, infinity, induction,\vspace{1mm}

\begin{tabular}{ll}
		\textbf{Separation} & $\exists y\, \forall z\, (z\in y\biimp z\in x\land \vp(z))$ 
		\ for all formulae $\vp$ \\
		\textbf{Collection} &  $\forall u\in x\, \exists v\, \vp(u,v)\imp \exists y\, \forall u\in x\, \exists v\in y\, \vp(u,v)$ \
		for all formulae $\vp$ \\
		\textbf{Powerset} &   $\exists y\, \forall z\, (\forall u\in z\, (u\in x)\imp z\in y)$
\end{tabular}

\medskip
	
Thus $\izf$ is  $\czf$ with bounded separation replaced by full separation and subset collection replaced by powerset. Note that powerset implies subset collection, and strong collection follows from separation and collection. It is well known that  $\czf$ with classical logic, and  hence $\izf$ with classical logic, is equivalent to $\sf{ZF}$.

	\subsection{Countable and dependent choice}
The full axiom of choice $\ac$  is known to imply forms of excluded middle \cite{D75}. For example, 
$\czf+\ac+\text{ separation}=\izf+\ac=\sf{ZFC}$ (cf.\ \cite{czf}). Nevertheless, countable choice and several forms of dependent choice  are deemed  constructively acceptable.  Note that $\rdc$ is an axiom of Myhill's $\sf CST$ (constructive set theory) \cite{M75}, which provides a standard set-theoretical framework for Bishop-style constructive mathematics \cite{B67}. Dependent choice also features in Brouwer's intuitionistic analysis (cf.\ \cite{TD88}). \vspace{1mm}

\textbf{Countable choice} ($\ac_{\omega}$): if $\forall n\in\omega\, \exists y\, \vp(n,y)$, then there exists a function $f$ with $\dom(f)=\omega$ such that $\forall n\in\omega\, \vp(n,f(n))$, for all formulae $\vp$.

\textbf{Dependent choice} ($\dc$):   if $\forall x\, \exists y\, \vp(x,y)$, then for every $x$ there is a function $f$ with $\dom(f)=\omega$ such that $f(0)=x$ and $\forall n\in\omega\, \vp(f(n),f(n+1))$, for all formulae $\vp$.

\textbf{Limited dependent  choice} ($\dc^\dagger$): if $X$ is a set and $\forall x\in X\, \exists y\in X\, \vp(x,y)$, then for every $x\in X$ there is a function $f\colon\omega\to X$ such that $f(0)=x$ and $\forall n\in\omega\, \vp(f(n),f(n+1))$, for all formulae $\vp$.

\textbf{Relativized dependent  choice} ($\rdc$): if $\forall x\, (\psi(x)\imp \exists y\, (\psi(y)\land \vp(x,y)))$, then for every $x$ such that $\psi(x)$, there is a function $f$ with $\dom(f)=\omega$ such that $f(0)=x$ and $\forall n\in\omega\, \vp(f(n),f(n+1))$, for all formulae $\psi$ and $\vp$. 

\textbf{Limited relativized dependent  choice} ($\rdc^\dagger$): if $X$ is a set and  $\forall x\in X\, (\psi(x)\imp \exists y\in X\, (\psi(y)\land \vp(x,y)))$, then for every $x\in X$ such that $\psi(x)$, there is a function $f\colon\omega\to X$ such that $f(0)=x$ and $\forall n\in\omega\, \vp(f(n),f(n+1))$, for all formulae $\psi$ and $\vp$. \vspace{1mm}

Note that classically, say over $\sf ZF$, the above forms of dependent choice are all equivalent to one another. Over $\czf$,  $\rdc\imp \dc$ and  $\rdc\imp\rdc^\dagger\imp \dc^\dagger\imp \ac_\omega$. Over $\izf$, due to separation, $\rdc^\dagger\biimp \dc^\dagger$.

\subsection{Inductive definitions and the regular extension axiom}
The theory $\czf$ allows for a smooth treatment of inductively defined classes   \cite{A86,czf}.  

\begin{definition}
	An inductive definition is a class $\Phi$ of ordered pairs. In other words,  any formula $\Phi(X,x)$ is an inductive definition.  We write \[ \rul{X}{x}{\Phi},\]
	or simply $(X,x)\in\Phi$, in place of $\Phi(X,x)$. A class $I$ is closed under $\Phi$ if, whenever $(X,x)\in\Phi$ and  $X\subseteq I$, we have $x\in I$. 
\end{definition}

\begin{theorem}
	Given an inductive definition $\Phi$, there exists a class $I(\Phi)$ such that 
	$\czf$ proves that $I(\Phi)$ is closed under $\Phi$, and for every class $J$, $\czf$ proves that if $J$ is closed under $\Phi$, then $I\subseteq J$. We say that $I(\Phi)$ is inductively defined by $\Phi$.
\end{theorem}
\begin{proof}
	See e.g.\  \cite{czf} for  details. The class $I(\Phi)$ is defined as follows. A set of ordered pairs $Z$ is called an iteration set for $\Phi$ if for every $\pair{u,x}\in Z$ there exists $(X,x)\in \Phi$ such that $X\subseteq \bigcup_{v\in u}\{z\mid \pair{v,z}\in Z\}$.	
	Let 
	\[ I(\Phi)=\bigcup_u\{ x\mid \text{there is an iteration set $Z$ such that }\pair{u,x}\in Z\}.\] 
	To prove  that $I(\Phi)$ is closed under $\Phi$ we invoke strong collection. The fact that $I(\Phi)$ is the smallest class closed under $\Phi$ is proved by induction.
\end{proof}
\begin{theorem}[Induction on the inductive definition  $\Phi$]
	$\czf$ proves
	\[ \forall (X,x)\in\Phi\, (\forall y\in X\, \vp(y)\imp \vp(x))\imp \forall x\in I(\Phi)\, \vp(x), \]
	for every formula $\vp$.
\end{theorem}
\begin{proof}
	Obvious.
\end{proof}

In the presence of Aczel's $\rea$ (regular extension axiom), one can show that many inductively defined classes exist as sets.

\begin{definition}[Regular]
	A set $X$ is {\em regular} if it is transitive and, whenever $x\in X$ and $R$ is a set such that  $\forall u\in x\, \exists v\in X\, (\pair{u,v}\in R)$, then there is a $y\in X$ such that
	\[ \forall u\in x\, \exists v\in y\, (\pair{u,v}\in R)\land \forall v\in y\, \exists u\in x\, (\pair{u,v}\in R). \]
	In particular, if $R\colon x\to X$ is a function, then $\ran(R)\in X$. 
\end{definition}

\textbf{Regular extension axiom} ($\rea$): every set is a subset of a regular set.	

\begin{definition}
	An inductive definition $\Phi$ is {\em bounded} if: 
	\begin{itemize}[leftmargin=5mm]
		\item the class $\{x\mid (X,x)\in\Phi\}$ is a set for every set $X$;
		\item there is a set $B$ such that whenever $(X,x)\in \Phi$, then $X$ is an image of a set in $B$, that is, there is a function $f$ from $Y$ onto $X$ for some $Y\in B$.  The set $B$ is called a bound for $\Phi$.
	\end{itemize}
\end{definition}
\begin{theorem}[$\czf+\rea$]\label{bounded}
	If $\Phi$ is a bounded inductive definition, then the class $I(\Phi)$ is a set. 
\end{theorem}
\begin{proof}
	See \cite{A86}[Theorem 5.2].
\end{proof}

\subsection{Dependent type constructions}

The usual type constructors $\Prod$ (dependent product), $\Sum$ (dependent sum), $\I$ (identity), and $\W$ have obvious set-theoretic analogues.

\begin{definition}[$\czf$]
	Let $F$ be a function   with $\dom(F)=X$.  Let
	\begin{align*} 
		\Prod(X,F)&= \prod_{x\in X}F(x)=\{f\colon X\to \bigcup_{x\in X}F(x)\mid \forall x\in X\, (f(x)\in F(x))\} && (\text{{\em cartesian product}}), \\
		\Sum(X,F)&=\sum_{x\in X}F(x)=\{ \pair{x,y}\mid x\in X\land y\in F(x)\} && (\text{{\em disjoint union}}).
	\end{align*}
	For $x$, $y$ sets, let $\I(x,y)=\{z\in \{0\}\mid x=y\}$.
\end{definition}
\begin{remark}
	Set exponentiation $X\to Y$ and the binary cartesian product $X\times Y$ are special instances of $\Prod(X,F)$ and $\Sum(X,F)$ respectively, with $F(x)=Y$.
\end{remark}

\begin{definition}[$\czf$]
	If  $F$ is a function with $\dom(F)=X$, then 
	$\W(X,F)$
	is the smallest class $W$ such that
	\begin{itemize}[leftmargin=5mm] 
		\item if $x\in X$ and $f\colon F(x)\to W$, then $\pair{x,f}\in W$.
	\end{itemize}
\end{definition}
\begin{remark}
	$\W(X,F)$ is the class of well-founded rooted trees with labels in $X$ such that the arity of a node with label $x$ is the set $F(x)$. 
\end{remark}

\begin{theorem}[$\czf+\rea$]
	The class $\W(X,F)$ is a set.
\end{theorem}
\begin{proof}
	By Theorem \ref{bounded}, since $\W(X,F)$ is bounded inductively definable.
\end{proof}

\begin{definition}[$\czf$]
	The  class $I=\Prod\Sum$ is inductively defined by clauses:
	\begin{enumerate}[label=(\roman*), leftmargin=*]		
		\item $n\in I$ for every $n\in\omega$;
		
		\item  $\omega\in I$;
		
		\item if $X\in I$ and $F\colon X\to I$, then  $\Prod(X,F),\Sum(X,F)\in I$.\vspace{1mm}
		
		The class $I=\Prod\Sum\I$ is inductively defined by clauses (i), (ii), (iii), and
		
		\item if $X\in I$  and $x,y\in X$, then $\I(x,y)\in I$.
	\end{enumerate}
\end{definition}

\begin{definition}[$\czf+\rea$]
	The class $I=\Prod\Sum\W$ is inductively defined by clauses (i), (ii), (iii) and
	\begin{enumerate}[label=(\roman*), leftmargin=*, start=5]
		\item if $X\in I$ and $F\colon X\to I$, then $\W(X,F)\in I$.
	\end{enumerate}
	
	The class $I=\Prod\Sum\W\I$ is inductively defined by clauses (i),(ii),(iii),(iv),(v).
\end{definition}

Given a class $I$, we say that $x$ is an $I$-set if $x$ belongs to $I$.

\subsection{Choice for dependent types}\label{sec:choice types}

\begin{definition}[Base]
	A set $X$  is a {\em base} if for every relation $R$ such that $\forall x\in X\, \exists y\, (\pair{x,y}\in R)$,  there exists a function $f$ with domain $X$ such that $\forall x\in X\, (\pair{x,f(x)}\in R)$.\footnote{In Aczel \cite{A86}, a set $X$ is said to be a base if  $\Prod(X,F)$ is inhabited, whenever   $F$  is a function with domain $X$ and $F(x)$ is inhabited for all $x\in X$, where a set $X$ is  inhabited if $\exists x\, (x\in X)$. Over $\czf$, the two definitions are equivalent.}
\end{definition}

\textbf{Presentation axiom} ($\pax$): every set is an image of a base.

$\Prod\Sum$ \textbf{axiom of choice} ($\Prod\Sum\dash\ac$): every $\Prod\Sum$-set is a base. 

$\Prod\Sum$ \textbf{presentation axiom} ($\Prod\Sum\dash\pax$): every set is an image of a $\Prod\Sum$-set and every $\Prod\Sum$-set is a base.\vspace{1mm}

We have analogue axioms for $\Prod\Sum\I$ and relatives.  A simple argument shows that $\czf+\pax\vdash \dc^\dagger$ (cf.\ \cite[p.\, 65]{A78} and \cite[p.\, 27]{A82}.)  Aczel's interpretation of $\czf$ in  type theory \cite{A78} validates $\rdc$ \cite[Theorem 7.1]{A78}\cite[Theorem 5.7]{A82}  and  $\Prod\Sum\I\dash\pax$ \cite[Theorem 7.4]{A82}. In \cite{A86},  the type-theoretic interpretation of $\czf+\rea$ (by using $\W$-types) is shown to validate $\rdc$ and $\Prod\Sum\W\I\dash\pax$  \cite[Theorem 5.6]{A86}.\footnote{The axiom $\rdc$ is simply  dubbed $\dc$ in \cite{A78,A82,A86}.}  Similarly, Rathjen's formulae-as-classes interpretation provides a model of $\czf+\rdc+\Prod\Sum\I\dash\ac$ in $\czf_{\sf exp}$ \cite[Theorem 4.13]{R03m} and of $\czf+\rea+\rdc+\Prod\Sum\W\I\dash\ac$ in $\czf_{\sf exp}+\rea$ \cite[Theorem 4.33]{R03m}, where $\czf_{\sf exp}$ is $\czf$ with exponentiation in lieu of subset collection. \vspace{1mm}

We introduce the following  natural extensions of $\ac_{\ft}$ and $\rdc_{\ft}$. \vspace{1mm}

$\Prod\Sum$ \textbf{limited axiom of choice} ($\Prod\Sum\dash\ac^\dagger$): if $X$ and $Y$ are $\Prod\Sum$-sets and $\forall x\in X\, \exists y\in Y\, \vp(x,y)$, then there is a function $f\colon X\to Y$ such that $\forall x\in X\, \vp(x,f(x))$, for all formulae $\vp$.

$\Prod\Sum$  \textbf{limited relativized dependent  choice} ($\Prod\Sum\dash\rdc^\dagger$):  this is just $\rdc^\dagger$ restricted to $\Prod\Sum$-sets.
\vspace{1mm}

We have similar axioms for $\Prod\Sum\I$ and relatives. Each schema  can be replaced by a single axiom over say $\izf+\rea$. This is just an upper bound. For example, over $\czf$, $\Prod\Sum\dash\ac^\dagger$ is equivalent to: \vspace{1mm} 

If  $X$ and $Y$ are $\Prod\Sum$-sets, then for every relation $R$ such that $\forall x\in X\, \exists y\in Y\, (\pair{x,y}\in R)$,  there is a function $f\colon X\to Y$ such that $\forall x\in X\, (\pair{x,f(x)}\in R)$.\vspace{1mm}

Clearly, $\Prod\Sum\dash\ac^\dagger$ follows from $\Prod\Sum\dash\ac$, and similarly for $\Prod\Sum\I$ and the like. Note that  $\Prod\Sum\dash\ac^\dagger$ already includes the axiom of choice for all finite types. 

\begin{theorem}\label{aczel}
	Over $\czf$, $\Prod\Sum\I\dash\ac\biimp \Prod\Sum\dash\ac$ and $\Prod\Sum\I\dash\pax$ $\biimp \Prod\Sum\dash\pax$.
	
	Over $\czf+\rea$, $\Prod\Sum\W\dash\ac\biimp \Prod\Sum\dash\ac$ and  $\Prod\Sum\W\dash\pax\biimp \Prod\Sum\dash\pax$
\end{theorem}
\begin{proof}
	See \cite[Theorem 3.7, Theorem 5.9]{A86}.
\end{proof}

Theorem \ref{aczel} does not  translate to $\Prod\Sum\dash\ac^\dagger$ and relatives. For example, the proof of  \cite[Theorem 3.7]{A86}  consists in showing, under $\Prod\Sum\dash\ac$, that every set in $\Prod\Sum\I$ is in bijection with a set in $\Prod\Sum$.  The use of $\Prod\Sum\dash\ac$ is essential. The proof of the following is trivial.

\begin{theorem}\label{image}
	Over $\czf +$\lq\lq every set is an image of a $\Prod\Sum\I$-set\rq\rq, 
	$ \Prod\Sum\I\dash\ac \biimp  \Prod\Sum\I\dash\ac^\dagger$.  
	Similarly for $\czf+\rea$ and $\Prod\Sum\W\I$ in place of $\czf$ and $\Prod\Sum\I$ respectively.
\end{theorem}

\section{Extensional generic realizability}\label{sec:realizability}

\subsection{Partial combinatory algebras}
Generic realizability is based on the notion of partial combinatory algebra. Let us review some basic facts. For more information, we refer the reader to \cite{B85,F75,F79,Oosten08}. 

\begin{definition}
	{\em Application terms}, or simply {\em terms}, over a set $A$ are defined by clauses:
	\begin{enumerate}[label=(\roman*), leftmargin=*]
		\item  variables $x_1,x_2,\ldots$  are terms;  
		\item   elements of $A$ are terms; 
		\item if $t$ and $s$ are terms,  $(ts)$ is also a term.
	\end{enumerate}
\end{definition}
\begin{definition}
	A {\em partial algebra} is a set $A$ together with  a partial binary operation $\cdot$ on $A$.
\end{definition}

Every partial algebra $A$ induces a partial interpretation of closed applications terms over $A$ by elements of $A$ in the obvious way. 

\begin{notation} 
	Let $a,b,c\in A$ and $t,s$ be closed terms over $A$. Instead of $a\cdot b$ we just write $ab$. We also employ the association to the left convention, meaning  that  $abc$ stands for $(ab)c$.  The relation   $t\downarrow a$ is defined by clauses: (i) $a\downarrow a$ and (ii) $(ts)\downarrow a$ if $t\downarrow b$, $s\downarrow c$ and $bc=a$ for some $b,c\in A$.    We write $t\downarrow$ for $\exists a\, (t\downarrow a)$.  We use Kleene equality $t\simeq s$ to express that either both sides are undefined, or else they are both defined and equal.  In the second case, we write $t=s$. 
\end{notation}

\begin{definition}
	A partial algebra $(A,\cdot)$ with at least two elements is a {\em partial  combinatory algebra}  (pca) if there are  elements 
	$\mb k$ and $\mb s$ in $A$ such that for all $a,b,c\in A$:
	\begin{itemize}[leftmargin=5mm]
		\item $\mb ka b\simeq a$;
		\item $\mb sab\downarrow$ and $\mb sabc\simeq ac(bc)$.
	\end{itemize}
\end{definition}
The combinators $\mb k$ and $\mb s$ are due   to Sch\"onfinkel \cite{Schoen24} and the defining equations, although formulated just in the total case,  are due to Curry \cite{Curry30}. The name combinatory stems from the property known as {\em combinatory completeness}. Informally, this means that we  can form   $\lambda$-terms.

\begin{lemma}[$\lambda$-abstraction]
	Let $A$ be a pca. For every term $t(x,x_1,\ldots,x_n)$, one can find (in an effective way) a new term $s(x_1,\ldots,x_n)$, denoted $\lambda x.t$, such that for all $a_1,\ldots,a_n\in A$:
	\begin{itemize}[leftmargin=5mm]
		\item $s(a_1,\ldots,a_n)\downarrow$;
		\item  $s(a_1,\ldots,a_n)a\simeq t(a,a_1,\ldots,a_n)$.
	\end{itemize}
\end{lemma}
\begin{remark}
	The term $\lambda x.t$ is built solely with the aid of $\mb k$, $\mb s$ and symbols occurring in $t$. Also, the construction of this term is uniform in the given pca.  
\end{remark}

An immediate consequence of $\lambda$-abstraction is the existence of pairing and unpairing combinators $\mb p$,  $\mb {p_0}$, and  $\mb {p_1}$ such that
$\mb pab\downarrow$ and  $\mb {p_i}(\mb pa_0a_1)\simeq a_i$.\footnote{Let $\mb p=\lambda xyz.zxy$, $\mb{p_0}=\lambda x.x\mb k$, and $\mb{p_1}=\lambda x.x\bar{\mb k}$, where $\bar{\mb k}=\lambda xy.y$. For our purposes, projections $\mb{p_0}$ and $\mb{p_1}$ need not be total.} A more remarkable  application of $\lambda$-abstraction is the recursion theorem for pca's.
\begin{lemma}[Recursion theorem] There exists a closed term $\mb f$ such that for all $a,b\in A$ we have $\mb f a\downarrow$ and  $\mb f ab\simeq a(\mb fa) b$. 
\end{lemma}

It is worth considering some additional structure (see  however Remark \ref{remark}).

\begin{definition}
	We say that $A$ is a pca over $\omega$ if there are combinators $\mb{succ}, \mb{pred}$ (successor and predecessor combinators), $\mb d$ (definition by cases combinator), and a  map
	$n\mapsto \bar n$ from $\omega$ to $A$ such that for all $n\in \omega$ 
	\begin{align*}
		\succe \bar n&\simeq \overline{n+1}, & \pred\overline{n+1}&\simeq \bar n, &
		\mb d\bar n\bar mab\simeq
		\begin{cases} a & n=m;\\ b & n\neq m. \end{cases}
	\end{align*}
	One then defines $\mb 0=_{\mathrm{def}}\bar 0$ and $\mb 1=_{\mathrm{def}}\bar 1$. 
\end{definition} 

For the rest of the paper, by a pca we really mean a pca over $\omega$. The only closed terms we ever need are built  up by using the combinators 
$\mb k$, $\mb s$, $\mb p$, $\mb{p_0}$, $\mb{p_1}$, $\mb 0$, $\mb 1$, $\mb{succ}$, $\mb{pred}$, $\mb d$. 

Note that one can do without $\mb k$ by letting $\mb k=_{\mathrm{def}}\mb d\mb 0\mb 0$. The existence of $\mb d$ implies that the map $n\mapsto \bar n$ is one-to-one. In fact, suppose $\bar n=\bar m$ but $n\neq m$. Then $\mb d\bar n\bar n\simeq \mb d\bar n\bar m$. It then follows that $a\simeq \mb d\bar n\bar nab\simeq \mb d\bar n\bar mab\simeq b$ for all $a,b$. On the other hand, by our definition, every pca contains at least two elements.

\begin{remark}\label{remark}
	Every pca can be turned into a pca over $\omega$.  For example, one can define Curry numerals and construct by $\lambda$-abstraction a combinator $\mb d$ for this representation of natural numbers. In practice, however, all natural examples of pca come already equipped with a canonical copy of the natural numbers and pertaining combinators.  
\end{remark}

\subsection{Defining realizability}
\begin{definition}[$\czf$]
	Given a pca $A$, the class $\va$  is inductively defined by the clause:
	\begin{itemize}[leftmargin=5mm]
		\item if $x\subseteq A\times A\times \va$, then $x\in\va$.
	\end{itemize} 
\end{definition}

The inductive definition of $\va$ within $\czf$ is on par with that of $\V(A)$  \cite[3.4]{R06}.

\begin{notation} We use $(a)_i$ or simply $a_i$ for $\mb {p_i}a$.   Whenever we write a term $t$, we assume that it is defined. In other words, a formula $\vp(t)$ stands for $\exists a\, (t\downarrow a\land \vp(a))$.
\end{notation}

\begin{definition}[Extensional realizability] We define the relation $a=b\fo \vp$, where $a,b\in A$ and $\vp$ is a realizability formula with parameters in $\va$. The atomic cases fall under the scope of definitions by transfinite recursion. 
	\begin{align*}   
		a=b&\fo x\in y && \Leftrightarrow && \exists z\, (\langle (a)_0,(b)_0,z\rangle \in y\land (a)_1=(b)_1\fo x=z)\\
		a=b& \fo x=y && \Leftrightarrow  &&\forall \langle c,d,z\rangle \in x\, ((ac)_0=(bd)_0\fo z\in y) \text{ and }  \\
		&&&&&  \forall \langle c,d,z\rangle \in y\, ((ac)_1=(bd)_1\fo z\in x)\\
		a=b& \fo \vp\land \psi && \Leftrightarrow && (a)_0=(b)_0\fo \vp \land (a)_1=(b)_1\fo \psi \\
		a=b& \fo \vp\lor\psi &&  \Leftrightarrow && (a)_0=(b)_0=\mb 0\land (a)_1=(b)_1\fo \vp \text{ or } \\ 
		&&&&&  (a)_0=(b)_0=\mb 1\land (a)_1=(b)_1\fo \psi \\
		a=b&\fo \neg\vp && \Leftrightarrow && \forall c, d\, \neg (c=d\fo \vp) \\
		a=b&\fo \vp\imp\psi && \Leftrightarrow && \forall c,d\,  (c=d\fo \vp\imp  ac=bd\fo \psi) \\
		a=b& \fo \forall x\in y\, \vp && \Leftrightarrow && \forall \langle c,d,x\rangle\in y\, (ac=bd\fo \vp) \\
		a=b&\fo \exists x\in y\, \vp && \Leftrightarrow && \exists x\, (\langle (a)_0,(b)_0,x\rangle \in y\land (a)_1=(b)_1\fo \vp)\\
		a=b& \fo \forall x\, \vp && \Leftrightarrow && \forall x\in \va\, (a=b\fo \vp) \\
		a=b& \fo \exists x\, \vp && \Leftrightarrow && \exists x\in \va\, (a=b\fo \vp)
	\end{align*}
\end{definition}

\begin{notation} We write $\fo \vp$ in place of $\exists a,b\in A\, (a=b\fo \vp)$ and  $a\fo \vp$ for $a=a\fo \vp$. 
\end{notation}

We will repeatedly use the following  internal pairing function $\vpair{x,y}$ in $\va$. 
\begin{definition}[Pairing]\label{pair}
	For $x,y\in\va$, let
	\begin{align*}
		\vset{x}&=\{\pair{\mb 0,\mb 0,x}\}, \\
		\vset{x,y}&=\{\pair{\mb 0,\mb 0,x},\pair{\mb 1,\mb 1,y}\}, \\
		\vpair{x,y}&=\{\pair{\mb 0,\mb 0,\vset{x}}, \pair{\mb 1,\mb 1,\vset{x,y}}\}. 
	\end{align*}
	
	Note that all these sets are in $\va$.
\end{definition}

We use  $\op(z,x,y)$ as abbreviation for the set-theoretic formula expressing  that $z$  is the ordered pair of $x$ and $y$. In standard  notation, $z=\pair{x,y}$. Ordered pairs can be defined as usual as $\pair{x,y}=_{\mathrm{def}}\{\{x\},\{x,y\}\}$.

\begin{lemma}\label{pairs}
	There  are closed terms ${\mb u}_i$ such that $\czf$ proves 
	\begin{align*}
		{\mb u}_0 & \fo \op(\vpair{x,y},x,y),\\
		{\mb u}_1 & \fo \vpair{x,y}=\vpair{u,v}\imp x=u\land y=v, \\
		{\mb u}_2 &\fo 	\op(z,x,y) \imp z=\vpair{x,y}.
	\end{align*}
\end{lemma}
\begin{proof}
	This is similar to \cite[3.2, 3.4]{M84}. 	
\end{proof}

\begin{theorem}[\cite{FR21}]\label{sound}
	For every theorem $\vp$ of $\czf+\ac_\ft$, there is a closed term $\mb t$ such that  $\czf$ proves $\mb t\fo \vp$.  Same for $\izf$.
\end{theorem}

\begin{remark}
	The realizability relation $\fo$ uses $A$ as a parameter. We should then write $\fo_A$.
	The soundness theorem claims that for every theorem $\vp$ of $T$ there is a closed term $\mb t$ such that $T\vdash \forall A\, (A\text{ is a pca}\imp \mb t\fo_A \vp)$. For special cases, where the pca $A$ is definable as a set in $T$ and $T\vdash A \text{ is a pca}$, we obtain $T\vdash \mb t\fo_A \vp$. For instance, in the case of Kleene's first algebra, we have that for every theorem $\vp$, there is a closed term $\mb t$ such that $T\vdash \exists n\in\omega\, (\mb t\downarrow n\land n\fo \vp)$. 
\end{remark}

\begin{definition} \label{notation}
	Write $T\models \vp$ for $T\vdash \exists a,b\in A\, (a=b\fo \vp)$. In general, we write $T\models S$ if $T\models \vp$ for every sentence $\vp$ of $S$.  We say that $T$ is {\em self-validating} if $T\models T$.
\end{definition}

\begin{theorem}\label{self}
	Let $T$ be any of the theories obtained by appending to either $\czf$ or $\izf$ some or all of the axioms $\rea$, $\ac_{\omega}$, $\dc$, $\dc^\dagger$, $\rdc$, $\rdc^\dagger$,  $\pax$. Then $T$ is self-validating.
\end{theorem}
\begin{proof}
	This holds in the case of generic realizability. See \cite[Theorem 6.2]{R06} for $\rea$ and  \cite[Theorem 10.1]{R06} for $\dc^\dagger$ (there called $\dc$), $\rdc$, and $\pax$. The other axioms are treated in a similar way.  In the case of extensional generic realizability, the proof is a nearly verbatim copy of that for generic realizability. 
	
	Let us sketch a proof for $\pax$. We wish to find a closed term $\mb e$ such that, over  $\czf+\pax$,  
	\[  \mb e\fo \forall y\, \exists x\, \exists f\, (x \text{ is a base and } f\colon x\twoheadrightarrow y). \]  
	By definition, this means that for every $y\in\va$ there are $x,f\in\va$ such that 
	\[   \mb e\fo x \text{ is a base and } f\colon x\twoheadrightarrow y. \]
	Let $y\in\va$. By $\pax$ in the background universe, there exists a base $B$ and a function $F$ from $B$ onto $y$. Say $F(u)=\pair{a_u,b_u,z_u}\in y$ for $u\in B$.  By transfinite recursion, define 
	\[  \check u=\{\pair{\mb 0,\mb 0, \check{v}}\mid v\in u\}. \]
	Then for every set $u$, $\check u\in\va$. Moreover, by induction, one can show that these names are absolute, namely,
	\[ u=v \biimp (\fo \check{u}=\check{v}). \]
	Let
	\begin{align*}
		x&=\{\pair{a_u,b_u,\check u}\mid u\in B\},\\
		f&=\{\pair{a_u,b_u,\vpair{\check u,z_u}}\mid u\in B\},
	\end{align*}
	where $\vpair{u,v}$ is the internal pair in $\va$ from Definition \ref{pair}.
	One can find a closed term $\mb f$  such  that $\mb f \fo \fun(f)$ thanks to the absoluteness of $\ \check{}\  $ names.  It is not difficult to build a closed term $\mb t$ such that $\mb t \fo f\colon x\twoheadrightarrow y$. Finally, it is straightforward to cook up a closed term $\mb b$ and show that $\mb b\fo x \text{ is a base}$,  by using the fact  that $B$ is a base.
\end{proof}

In view of Definition \ref{notation}, our goal is to show that
\begin{itemize}[leftmargin=5mm]
	\item $T\models T+\Prod\Sum\I\dash\ac^\dagger+\Prod\Sum\I\dash\rdc^\dagger$,
	\item $T+\rea\models T+\rea+\Prod\Sum\W\I\dash\ac^\dagger+\Prod\Sum\W\I\dash\rdc^\dagger$,
\end{itemize}
where $T$ is as in Theorem \ref{self}.

\section{Realizing choice for dependent types}\label{sec:choice}

\subsection{Dependent types over partial combinatory algebras}\label{sec:types}

Given a pca $A$, we wish to define extensional type structures on $A$, that is,  structures of the form 
\[   (\mathbb{T},\sim,(A_\sigma,\sim_\sigma)_{\sigma\in\mathbb{T}}), \]
where $\sim$ is a  partial equivalence relation   on $A$  and $\sim_\sigma$  is a  partial equivalence relation on $A$ for every $\sigma$ such that $\sigma\sim\sigma$. Then $\mathbb{T}=\{\sigma\in A\mid \sigma\sim\sigma\}$ and $A_\sigma=\{a\in A\mid a\sim_\sigma a\}$  for every  $\sigma\in\mathbb{T}$. Clearly, we want equivalent types to have the same elements. We must then ensure that  if $\sigma\sim\tau$,  then  $a\sim_\sigma b$ iff $a\sim_\tau b$.

\begin{definition}
	We say that $\sigma$ is a type if $\sigma\in\mathbb{T}$. We say that $a\in A$ has type $\sigma$ if $a\in A_\sigma$.
\end{definition}

From now on, definitions and theorems take place in $\czf$, unless we are dealing with $\W$-types, in which case we work in $\czf+\rea$. 

\begin{notation}
	Given a pca $A$, $a,b,c\in A$, and $n\in\omega$, we write:
	\begin{align*}
		{\sf N}_n&=_{\mathrm{def}} \mb p \bar 0\bar n, & {\sf N}&=_{\mathrm{def}}\mb p\bar 1\mb 0, \\
		\Prod_ab&=_{\mathrm{def}} \mb p \bar 2(\mb pab), &  \Sum_ab&=_{\mathrm{def}}\mb p \bar 3(\mb pab), \\
		\I_c(a,b)&=_{\mathrm{def}} \mb p \bar 4(\mb p c(\mb pab)),  &
		\W_ab &=_{\mathrm{def}}\mb p \bar 5(\mb p ab).
	\end{align*}
	This  ensures unique readability. For example, ${\sf N}\neq \Prod_a b$, and $\Prod_ab=\Prod_cd$ implies $a=b$ and $c=d$. 
\end{notation}

\begin{definition}[$\Prod\Sum\W\I$-types]\label{types}
	Given a pca $A$, we inductively define the $\Prod\Sum$ type structure   by clauses:
	\begin{enumerate}[label=(\roman*), leftmargin=*]
		\item ${\sf N}_n\sim {\sf N}_n$ and $a\sim_{{\sf N}_n} b$ iff $a=b=\bar m$ for some $m<n$;
		
		\item ${\sf N}\sim {\sf N}$ and $a\sim_{\sf N} b$ iff $a=b=\bar n$ for some $n\in\omega$;
		
		\item if $\sigma\sim\tau$, and $a\sim_\sigma b$ implies   $ia\sim jb$, then
		\begin{itemize}[leftmargin=5mm]
			\item  $\alpha=\Prod_\sigma i\sim \Prod_\tau j$;
			\item  $f\sim_\alpha g$ iff $fa\sim_{ia} gb$ whenever $a\sim_\sigma b$;
			\item[] and 
			\item $\beta=\Sum_\sigma i\sim \Sum_\tau j$; 
			\item $a\sim_\beta b$ iff $a_0\sim_\sigma b_0$ and $a_1\sim_{ia_0} b_1$.
		\end{itemize}
	\end{enumerate}
	
	The $\Pi\Sum\I$ type structure   is inductively defined by adding the clause:
	\begin{enumerate}[label=(\roman*), leftmargin=*, start=4]	
		\item if $\sigma\sim\tau$, $a\sim_\sigma \breve a$, and $b\sim_\sigma \breve b$, then
		\begin{itemize}[leftmargin=5mm]
			\item $\gamma=\I_\sigma(a,b)\sim \I_\tau(\breve a,\breve b)$;
			\item $c\sim_\gamma d$ iff $c=d=\mb 0$ and $a\sim_\sigma b$.
		\end{itemize}
	\end{enumerate}
	
	The $\Prod\Sum\W\I$  type structure  is inductively defined by clauses (i), (ii), (iii), (iv), and 
	\begin{enumerate}[label=(\roman*), leftmargin=*, start=5]	
		\item if $\sigma\sim\tau$, and $a\sim_\sigma b$ implies   $ia\sim jb$, then
		\begin{itemize}[leftmargin=5mm]
			\item  $\delta=\W_\sigma i\sim\W_\tau j$;
			\item $c\sim_\delta d$ iff $c_0\sim_\sigma d_0$ and $c_1p\sim_\delta d_1q$ whenever  $p\sim_{ic_0} q$.
		\end{itemize}
	\end{enumerate}
\end{definition}

\begin{remark}\label{set} 
	To be precise, each type structure from Definition \ref{types} is an inductively defined  class $I$  of triples $\pair{\sigma,\tau,B}$. The intended meaning is $\sigma\sim\tau$ and $a\sim_\sigma b$ iff $\pair{a,b}\in B$. 
	For instance, the rules for the introduction of   $\Prod$-types  and $\W$-types are 
	\begin{equation*} 
		\rul{X}{\pair{\Prod_\sigma i,\Prod_\tau j,C}}{}\  \text{ and }\  \rul{X}{\pair{\W_\sigma i,\W_\tau j, D}}{},  
	\end{equation*}
	where $\sigma,\tau,i,j\in A$, and for some  set $B\subseteq A\times A$ and some function  $F$ with $\dom(F)=B$, we have
	\begin{itemize}[leftmargin=5mm]
		\item  $\{\pair{\sigma,\tau,B}\}\cup\{\pair{ia,jb,F(a,b)}\mid \pair{a,b}\in B\}\subseteq X$,  
		\item $C=\{\pair{f,g}\in A\times A\mid \forall \pair{a,b}\in B\, (\pair{fa,gb}\in F(a,b))\}$, 
		\item  $D=\{\pair{c,d}\in A\times A\mid \pair{c_0,d_0}\in B\land \forall \pair{p,q}\in F(c_0,d_0)\, \pair{c_1p,d_1p}\in D\}$.
	\end{itemize}
	It is not difficult to see that $D$ has a bounded inductive definition, and therefore by Theorem \ref{bounded} is a set in $\czf+\rea$.
	
	By induction on the inductive definition of $I$, one may show that if $\pair{\sigma,\tau_0,B_0}\in I$ and $\pair{\sigma,\tau_1,B_1}\in I$, then $B_0=B_1$. Hence, by letting  
	\begin{align*}
		\sim &=_{\mathrm{def}} \{\pair{\sigma,\tau}\mid \exists B\, \pair{\sigma,\tau,B}\in I\}, \\
		\sim_\sigma &=_{\mathrm{def}} \bigcup\{ B\mid \exists \tau\, \pair{\sigma,\tau,B}\in I\},
	\end{align*}
	we obtain in particular that if $\pair{\sigma,\tau,B}\in I$, then $B=\bigcup \{B\mid \exists \tau\, \pair{\sigma,\tau,B}\in I\}$ is a set, and so $\sim_\sigma$ is a set as well.
\end{remark}

\begin{lemma}\label{per}
	Suppose $\sigma\sim\tau$ and $\tau\sim\rho$. Then
	\begin{enumerate}[label=\textup{(\roman*)}, leftmargin=*]
		\item  $\sigma\sim\sigma$, $\tau\sim\sigma$ and $\sigma\sim\rho$;
		\item  $a\sim_\sigma b$ iff $a\sim_\tau b$;
		\item $a\sim_\sigma b$ and $b\sim_\sigma c$ implies $a\sim_\sigma a$, $b\sim_\sigma a$, and $a\sim_\sigma c$.
	\end{enumerate}
	In other words, $\sim$ and $\sim_\sigma$ are partial equivalence relations on $A$.
\end{lemma}
\begin{proof}
	By induction on the (inductive definition of the) type structure. The case of $\W_\sigma i$ is dealt with a further induction on the given $\W$-type. 
\end{proof}

\subsection{Representing $\Prod\Sum\W\I$}

\begin{definition}
	Given a type structure on $A$ and a type $\sigma$, we say that $i\in A$ is a {\em family of types over} $\sigma$ if 	 $ia\sim ib$ whenever $a\sim_\sigma b$.  If $\sigma$ is a type and $ia\sim jb$ whenever $a\sim_\sigma b$, we write $i\approx_\sigma j$.
\end{definition}

Note that in any of the type structures here considered,  if $i$ is a family of types over $\sigma$,  then $\Prod_\sigma i$ and $\Sum_\sigma i$ are types. Moreover,  if $\sigma\sim\tau$ and $i\approx_\sigma j$, then $i\approx_\tau j$. 

\begin{definition}[Mapping typed elements of $A$ in $\va$]
	For $n\in\omega$, let 
	\[  \dot n=\{\pair{\bar m,\bar m,\dot m}\mid m<n\}. \]
	
	For every  $\Prod\Sum\I$-type $\sigma$, we inductively define $a^\sigma\in \va$ for every $a$ of type $\sigma$ by:
	\begin{itemize}[leftmargin=5mm]
		\item if $\sigma={\sf N}_n$ or $\sigma={\sf N}$, let $a^\sigma=\dot m$, where $a=\bar m$;
		\item if $\alpha=\Prod_\sigma i$, let 
		$   f^\alpha=\{\pair{a,b,\vpair{a^\sigma,(fa)^\tau}}\mid a\sim_\sigma b\land ia=\tau\}$; 
		\item if $\beta=\Sum_\sigma i$, let
		$   a^\beta=\vpair{(a_0)^\sigma,(a_1)^\tau}$,  where $ia_0= \tau$; 
		\item if $\gamma=\I_\sigma(a,b)$, let 
		$       c^\gamma=0 $.
	\end{itemize}
	
	For  $\Prod\Sum\W\I$-types, we add the clause:
	\begin{itemize}[leftmargin=5mm]
		\item if $\delta=\W_\sigma i$, let
		$c^\delta=\vpair{a^\sigma,\{\pair{p,q,\vpair{p^{ia},(c_1p)^\delta}}\mid p\sim_{ia} q\}}$, where $a=c_0$.
	\end{itemize}
\end{definition} 

\begin{remark}
	A comment on the above definition may be in order. The class function $(\sigma,a)\mapsto a^\sigma$ is obtained by  defining a  class $E$ of pairs  $\pair{\sigma,e_\sigma}$, the intended meaning of which is $e_\sigma=\{\pair{a,a^\sigma}\mid a\in A_\sigma\}$, where $A_\sigma=\{a\in A\mid a\text{ has type } \sigma\}$. Note that $A_\sigma$ is a set.  Hence, $a^\sigma$ is really $e_\sigma(a)$. For example, the rules for $\Sum$-types and $\W$-types are 
	\[ \rul{X}{\pair{\Sum_\sigma i,e_\beta}}{} \  \text{ and }\   \rul{X}{\pair{\W_\sigma i,e_\delta}}{}, \]
	where $\sigma$ is type, $i$ is a family of types over $\sigma$, and there exist functions  $e_\sigma$ and $f$ on $A_\sigma$ such that
	\begin{itemize}[leftmargin=5mm]
		\item  $\fun(f(a))\land \dom(f(a))=A_{ia}$, for every $a\in A_\sigma$,
		\item $X=\{\pair{\sigma,e_\sigma}\}\cup\{\pair{ia,f(a)} \mid a\in A_\sigma\}$,   
		\item $e_\beta=\{\pair{a,\vpair{e_\sigma(a_0),f(a_0,a_1)}}\mid a\text{ has type } \Sum_\sigma i\}$, 
	\end{itemize}
	and  $e_\delta$ is inductively defined by the clause:
	\begin{itemize}[leftmargin=5mm]
		\item if $c$ has type $\W_\sigma i$ and there exists a function $k$ with $\dom(k)=A_{ic_0}$ such that $\pair{c_1p,k(p)}\in e_\delta$ for every $p\in\dom(k)$, then 
		\[ \pair{c,\vpair{e_\sigma(c_0),\{\pair{p,q,\vpair{f(c_0,p), k(p)}}\mid p\sim_{ic_0} q\}}} \in e_\delta. \]
	\end{itemize} 
	By induction on the inductive definition of $E$, one shows that $E$ is a class function on $\mathbb{T}$, and $E(\sigma)=e_\sigma\colon A_\sigma\to \va$, for every type $\sigma$.  
	Note that $e_\delta$ has a bounded inductive definition and so is a set in $\czf+\rea$. 	
\end{remark}

\begin{definition}[Canonical $\Prod\Sum\W\I$-sets in $\va$]\label{nice}
	Given a  type $\sigma$ and a family $i$ of types over $\sigma$,   we  define sets $X_\sigma, F_{\sigma,i}\in\va$ by:
	\begin{itemize}[leftmargin=5mm]
		\item $X_\sigma=\{\pair{a,b,a^\sigma}\mid a \sim_\sigma b\}$;
		\item $F_{\sigma,i}=\{\pair{a,b,\vpair{a^\sigma,X_\tau}}\mid a\sim_\sigma b\land ia= \tau\}$.
	\end{itemize}
\end{definition}

\begin{remark}
	The set $X_{{\sf N}_n}=\dot n$ is the canonical name for $n\in\omega$, and
	$X_{\sf N}=\dot \omega=\{\pair{\bar n,\bar n, \dot n}\mid n\in\omega\}$ is the canonical name for $\omega$. 
\end{remark}

\begin{lemma}[Absoluteness and uniqueness up to extensional equality]\label{inj}
Let $\sigma\sim\tau$ and $a,b$ of type $\sigma$. Then:
\begin{itemize}[leftmargin=5mm]
	\item   $\fo  a^\sigma=b^\tau$ implies $a\sim_\sigma b$;
	\item $a\sim_\sigma b$ implies $a^\sigma=b^\tau$.
\end{itemize}
\end{lemma}
\begin{proof}
By induction on the (inductive definition of the) type structure. The cases $\sigma=\tau={\sf N}_n$ and $\sigma=\tau={\sf N}$ are straightforward. 

Let $\alpha=\Prod_\sigma i\sim \Prod_\tau j=\gamma$ and $f,g\in A_\alpha$. 
Recall that:
\begin{align*}
	f^\alpha&=\{\pair{a,b,\vpair{a^\sigma,(fa)^{ia}}}\mid a\sim_\sigma b\};\\
	g^\gamma&=\{\pair{a,b,\vpair{a^\tau,(ga)^{ja}}}\mid a\sim_\tau b\}. 
\end{align*}

Suppose
$\fo f^\alpha=g^\gamma$. 
Let us show that $f\sim_\alpha g$. Let $a\sim_\sigma b$. We want to prove  $fa\sim_\rho gb$, where $\rho=ia$.
By definition of realizability, we have that
\[ \fo \vpair{a^\sigma,(fa)^\rho}=\vpair{\breve{a}^\tau,(g\breve{a})^\eta}, \]
for some $\breve{a}\in A_\tau$, where $\eta=j\breve{a}$. Then
$\fo a^\sigma=\breve a^\tau$ and $\fo (fa)^\rho=(g\breve a)^\eta$.
By induction, $a\sim_\sigma \breve a$ and $fa\sim_\rho g\breve a$. By Lemma \ref{per}, $\breve a\sim_\sigma b$. Hence, $\breve a\sim_\tau b$, and therefore $g\breve a\sim_\eta gb$. Since $i\approx_\sigma j$, it also follows from $a\sim_\sigma \breve a$ that $\rho\sim \eta$. Thus $g\breve a\sim_\rho gb$. By transitivity, $fa\sim_\rho gb$. 

Now suppose $f\sim_\alpha g$. We wish to show that $f^\alpha=g^\gamma$. Since $\sigma\sim\tau$, it suffices to show that if $a\in A_\sigma$, then $a^\sigma=a^\tau$ and $(fa)^{ia}=(ga)^{ja}$. This follows by induction since $\sigma\sim\tau$, $a\sim_\sigma a$, $ia\sim ja$ and  $fa\sim_{ia} ga$.

The cases $\Sum_\sigma i$, $\I_\sigma(a,b)$ and $\W_\sigma i$ are left as an exercise.  
\end{proof}

By virtue of Lemma \ref{per} and Lemma \ref{inj}, every $X_\sigma$ is injectively presented in the following sense.

\begin{definition} 
A set $x\in\va$ is {\em injectively presented} if
\[   x=\{\pair{a,b,f(a)}\mid a\sim b\}, \]
where $\sim$ is a partial equivalence relation on $A$ and $f$ is a (definable) function  from the domain of $\sim$ to $\va$ such that:
\begin{itemize}[leftmargin=5mm]
	\item $\fo f(a)=f(b)$ implies $a\sim b$;
	\item $a\sim b$ implies $f(a)=f(b)$.
\end{itemize}
\end{definition}
This property will ensure the validity of  choice for names of the form $X_\sigma$ (see Theorem \ref{canonical choice}).

\begin{lemma}\label{inv}
Let $\sigma\sim\tau$ and $i\approx_\sigma j$. Then
$X_\sigma=X_\tau$ and
$F_{\sigma,i}=F_{\tau,j}$.
\end{lemma}
\begin{proof}
By Lemma \ref{inj}.
\end{proof}

\begin{theorem}\label{realizers}
There are closed terms such that:
\begin{align*}
	\tag{i} \mb o&\fo X_{\sf N}=\omega; \\
	\tag{ii} \mb{fun}&\fo \fun(F_{\sigma,i})\land \dom(F_{\sigma,i}) = X_\sigma;\\
	\tag{iii}      \mb{prod}&\fo X_{\Prod_\sigma i}=\Prod(X_\sigma, F_{\sigma,i}); \\
	\tag{iv}          \mb{sum}&\fo X_{\Sum_\sigma i}=\Sum(X_\sigma, F_{\sigma,i});  \\
	\tag{v}           \mb{id}&\fo X_{\I_\sigma(a,b)}=\I(a^\sigma,b^\sigma);  \\
	\tag{vi}          \mb w&\fo X_{\W_\sigma i}=\W(X_\sigma,F_{\sigma,i}),
\end{align*}
where $\sigma$ is a type, $i$ is a family of types over $\sigma$, and $a,b$ have type $\sigma$.	
\end{theorem}
\begin{remark}
$\czf$ proves (i), (ii), (iii), (iv), (v)  for $\Prod\Sum\I$-types. $\czf+\rea$ proves (i), (ii), (iii), (iv), (v), (vi) for $\Prod\Sum\W\I$-types. 
\end{remark}

\begin{theorem}\label{canonical}
There are closed terms $\mb i$ and $\mb e$ such that:
\[ c=d\fo X\in \Prod\Sum\I \ \ \text{implies} \ \ \sigma=\mb ic\sim\mb id\ \ \text{and}\ \  \mb ec=\mb e d\fo X=X_\sigma, \]
for every $X\in\va$. Same for $\Prod\Sum\W\I$.
\end{theorem}

We prove Theorem \ref{realizers} and  Theorem \ref{canonical} in the appendix.

\subsection{Realizing choice}

\begin{theorem}[$\ac_{\sigma,\tau}$]\label{canonical choice}
There is a closed term $\mb e$ such that $\czf$ proves
\[  \mb e \fo \forall x\in X_\sigma\, \exists y\in X_\tau\, \vp(x,y)\imp \exists f\colon X_\sigma\to X_\tau\, \forall x\in X_\sigma\, \vp(x,f(x)), \]
for all $\Prod\Sum\I$-types $\sigma,\tau$ and for every formula $\vp$. Same for $\czf+\rea$ and $\Prod\Sum\W\I$-types.
\end{theorem}
\begin{proof}
Suppose $c=d\fo \forall x\in X_\sigma\, \exists y\in X_\tau\, \vp(x,y)$. Then, whenever $a\sim_\sigma b$, we have  $(ca)_0\sim_\tau (db)_0$ and  	
$(ca)_1=(db)_1\fo \vp(a^\sigma,e^\tau)$,
where $e=(ca)_0$. Let
\[ f=\{\pair{a,b,\vpair{a^\sigma,e^\tau}}\mid a\sim_\sigma b\land e=(ca)_0\}. \]
Then $f\in\va$. We wish to find  $\mb e$ such that
\[  \mb e c=\mb ed\fo f\colon X_\sigma\to X_\tau\land  \forall x\in X_\sigma\, \vp(x,f(x)). \]
It suffices to look for closed terms $\mb{r},\mb{c},\mb{f}$ such that 
\begin{align*}
	\tag{1}	\mb{r}c=\mb{r}d & \fo f\subseteq X_\sigma\times X_\tau, \\
	\tag{2}   \mb cc=\mb cd &\fo \forall x\in X_\sigma\, \exists y\, \exists z\in f\, (\op(z,x,y)\land \vp(x,y)), \\
	\tag{3}   \mb fc=\mb fd &\fo f \text{ is functional},
\end{align*}
where $f\subseteq X_\sigma\times X_\tau$ stands for $\forall z\in f\, \exists x\in X_\sigma\, \exists y\in X_\tau\, \op(z,x,y)$, and $f$ is functional stands for
\[ \forall z_0\in f\, \forall z_1\in f\, \forall x\, \forall y_0\, \forall y_1\, (\op(z_0,x,y_0)\land \op(z_1,x,y_1)\imp y_0=y_1). \]

(1) and (2) are straightforward.  (3) Let $a\sim_\sigma b$ and $\breve a\sim_\sigma \breve b$. We want 
\[  \mb fca\breve a=\mb fd b\breve b\fo \op(\vpair{a^\sigma,e^\tau},x,y_0)\land \op(\vpair{\breve a^\sigma,\breve e^\tau},x,y_1)\imp y_0=y_1, \]
for all $x,y_0,y_1\in\va$, where $e=(ca)_0$ and $\breve e=(c\breve a)_0$.
Suppose 
\[  \breve c=\breve d\fo \op(\vpair{a^\sigma,e^\tau},x,y_0)\land \op(\vpair{\breve a^\sigma,\breve e^\tau},x,y_1). \]
In particular, $\fo a^\sigma=\breve a^\sigma$. By Lemma \ref{inj}, $a\sim_\sigma \breve a$. By Lemma \ref{per}, since $(ca)_0\sim_\tau (da)_0$ and $(c\breve a)_0\sim_\tau (da)_0$, it follows that $(ca)_0\sim_\tau (c\breve a)_0$. That is, $e\sim_\tau \breve e$. By Lemma \ref{inj} again, we get $e^\tau=\breve e^\tau$. By the properties of equality, there is a closed term $\mb q$ such that
$\mb q\breve c=\mb q\breve d\fo y_0=y_1$.
We  may let  $\mb f=_{\mathrm{def}}\lambda ca\breve a\breve c. \mb q\breve c$.
\end{proof}

\begin{theorem}[$\rdc_{\sigma}$]\label{canonical dependent}
There is a closed term $\mb e$ such that $\czf$ proves
\begin{multline*}  \mb e \fo \forall x\in X_\sigma\, (\psi(x)\imp \exists y\in X_\sigma\, (\psi(x)\land \vp(x,y))) \imp \forall x\in X_\sigma\, (\psi(x)\imp \\ \exists f\colon \omega\to X_\sigma\, (f(0)=x\land \forall n\in\omega\, \vp(f(n),f(n+1)))), \end{multline*}
for every $\Prod\Sum\I$-type $\sigma$ and for every formula $\vp$. Same for $\czf+\rea$ and $\Prod\Sum\W\I$-types.
\end{theorem}
\begin{proof}
Suppose 
\begin{equation}
	c=d\fo \forall x\in X_\sigma\, (\psi(x)\imp \exists y\in X_\sigma\, (\psi(y)\land \vp(x,y))), \tag{1}	
\end{equation}   
$\mathring a\sim_\sigma \mathring b$ and $\mathring c=\mathring d\fo \psi(a^\sigma)$. We search for $\mb e$ such that
\[  \mb e c\mathring a\mathring c=\mb ed\mathring b\mathring d\fo \exists f \colon \omega\to X_\sigma\, (f(0)=\mathring a^\sigma\land \forall n\in\omega\, \vp(f(n),f(n+1))). \]
By (1), whenever $a\sim_\sigma b$ and $\breve c=\breve d\fo \psi(a^\sigma)$, we have 
\begin{equation}   
	e=(ca\breve c)_0\sim_\sigma (d b\breve d)_0, \tag{2}
\end{equation}
\begin{equation}
	(ca\breve c)_1= (db\breve d)_1\fo \psi(e^\sigma)\land  \vp(a^\sigma, e^\sigma). \tag{3} 
\end{equation}
In every pca one can simulate primitive recursion. We thus have a closed term $\mb r$ satisfying the following equations, where  $r$ is an abbreviation for $\mb rc\mathring a\mathring c$:
\begin{align*}
	r\mb 0&\simeq \mb p \mathring a(\mb p\mathring c\mb 0);\\
	r\overline{n+1}&\simeq \mb p (c(r\bar n)_0(r\bar n)_{10} )_0  (c(r\bar n)_0(r\bar n)_{10})_1. 
\end{align*}
By induction, owing to (2) and (3), one can verify that for every $n\in\omega$,  
\[  a_n=(\mb rc\mathring a\mathring c \bar n)_0\sim_\sigma (\mb rd\mathring b\mathring d \bar n)_0=b_n, \]
\[  (\mb rc\mathring{a}\mathring{c}\: \overline{n+1})_1= (\mb rd\mathring b\mathring{d}\: \overline{n+1})_1\fo \psi(a_n^\sigma)\land \vp(a_n^\sigma,a_{n+1}^\sigma). \]
Note that 
$\mathring a=(\mb rc\mathring a\mathring c \mb 0)_0\sim_\sigma (\mb rd\mathring b\mathring d \mb 0)_0=\mathring b$, and 
$\mathring c=(\mb rc\mathring a\mathring c \mb 0)_{10}= (\mb rd\mathring b\mathring d \bar 0)_{10}=\mathring d \fo \psi(\mathring a^\sigma)$.
Let 
\[ f=\{\pair{\bar n,\bar n,\vpair{\dot n,a_n^\sigma}}\mid n\in\omega\land a_n=(\mb rc\mathring a\mathring c \bar n)_0\}. \]
Then $f\in\va$. By using $\mb r$, it is not difficult to build $\mb e$ such that 
\[  \mb e c\mathring a\mathring c=\mb ed\mathring b\mathring d\fo f \colon \omega\to X_\sigma\land  f(0)=\mathring a^\sigma\land \forall n\in\omega\, \vp(f(n),f(n+1)). \]
\end{proof}

\begin{theorem}\label{main}
There is a closed term $\mb c$ such that $\czf$ proves $\mb c\fo \chi$, for every instance $\chi$ of $\Prod\Sum\I\dash\ac^\dagger$. Same for  $\Prod\Sum\I\dash\rdc^\dagger$.  Similarly, there is a closed term $\mb c$ such that   $\czf+\rea$ proves $\mb c\fo \chi$, for every instance $\chi$ of  $\Prod\Sum\W\I\dash\ac^\dagger$. Same for $\Prod\Sum\W\I\dash\rdc^\dagger$. 
\end{theorem}
\begin{proof}
By Theorem \ref{canonical}, Theorem \ref{canonical choice} and Theorem \ref{canonical dependent}.
\end{proof}

\begin{corollary} It  thus follows by Theorem \ref{self} that 
\begin{enumerate}[label=\textup{(\roman*)}, leftmargin=*]
	\item $T\models T+\Prod\Sum\I\dash\ac^\dagger+\Prod\Sum\I\dash\rdc^\dagger$,
	\item  $T+\rea\models T+\rea+\Prod\Sum\W\I\dash\ac^\dagger+\Prod\Sum\I\dash\rdc^\dagger$,
\end{enumerate}
for every theory $T$ obtained by appending to either $\czf$ or $\izf$ some or all of the axioms $\ac_{\omega}$, $\dc$, $\dc^\dagger$, $\rdc$, $\rdc^\dagger$, $\pax$. 
\end{corollary}

\subsection{Failure of presentation}\label{sec: remarks}
The interpretations of $\czf$ and $\czf+\rea$ in type theory validate strong presentation axioms (see Section \ref{sec:choice types}). It is natural to ask whether  extensional generic realizability validates the statement that every set is an image of a $\Prod\Sum\I$-set, even under the assumption that this holds in the background universe.  If it did, in view of Theorems \ref{image} and \ref{main}, we would obtain that
$T+\text{every set is an image of a $\Prod\Sum\I$-set} \models \Prod\Sum\I\dash\pax$,
and similarly for $T+\rea$ and $\Prod\Sum\W\I$-sets. However, the answer is negative.

\begin{theorem} The following holds.
\[ \czf\models\text{there is a set that is not an image of a $\Prod\Sum\I$-set}. \]
Similarly for $\czf+\rea$ and $\Prod\Sum\W\I$-sets.
\end{theorem}
\begin{proof}
Let $y=\{\pair{\mb 0,\mb 0,\dot 0},\pair{\mb 0,\mb 0,\dot 1}\}\in\va$. Although $y$ is in bijection with $2=\{0,1\}$, this fact is not realizable.  The proof is the same in both cases. Let us consider the $\Prod\Sum\I$ case and show that
\[   \mb 0\fo \neg \exists X\, \exists f\,  (X\in\Prod\Sum\I\land f\colon X \twoheadrightarrow y). \]
Towards a contradiction, suppose $\fo X\in \Prod\Sum\I\land f\colon X \twoheadrightarrow y$, for some $X,f\in\va$. By soundness and Theorem \ref{canonical}, there is a $\Prod\Sum\I$-type $\sigma$ such that $\fo f\colon X_\sigma \twoheadrightarrow y$. In particular, there are $a,b$ such that
\[  a=b\fo \forall v\in y\, \exists x\in X_\sigma\, \exists z\in f\, \op(z,x,v). \]
By unpacking the definition, we must have  $\breve a= (a\mb 0)_0\sim_\sigma(b\mb 0)_0=\breve b$ and
\[ \fo \exists z\in f\, \op(z,x,\dot 0)\land \exists  z\in f\, \op(z,x,\dot 1), \]
where $x=\breve a^\sigma$. From $\fo \fun(f)$, we then obtain $\fo \dot 0=\dot 1$, which implies $0=1$, for the desired contradiction.
\end{proof}

\section{Far beyond Goodman's theorem}\label{sec:goodman}

Goodman \cite{G76,G78} proved that  intuitionistic finite type arithmetic ${\sf HA}^\omega$ augmented with the axiom of choice in all finite types $\ac_{\ft}$  is conservative over $\sf HA$. This still holds in the presence of extensionality,  that is,  $\sf{HA}^\omega+{\sf EXT}+{\sf AC_{\ft}}$ is conservative over $\sf HA$ (cf.\ Beeson \cite{Beeson79,B85}). Since ${\sf HA}^\omega$ and ${\sf HA}^\omega +{\sf EXT}$ are easily seen to be conservative over $\sf HA$ (via the models $\sf HRO$ and $\sf HEO$ respectively, cf.\ \cite{T73}), Goodman's theorem is the prototypical example of arithmetical conservativity  of (some) choice over a given theory $T$. Goodman-type results have been proved for  a variety of set theories and choice principles at finite types. Write $T=_{arith}S$ to convey that $T$ and $S$ prove the same arithmetic sentences.
\begin{theorem} This much (as a minimum) is known:
\begin{itemize}[leftmargin=5mm]
	\item $T =_{arith}T+\cac_{\ft}+{\dc}_{\ft}$ for $T=\izf$ and extensions of $\izf$ by large set axioms, see Friedman and  Scedrov \cite[Theorem 5.1]{friedman_scedrov84};
	\item $T =_{arith}T+\ac_{\ft}$ for $T=\czf$ and fragments thereof, see Gordeev \cite[p.\ 25]{gordeev}.
\end{itemize}
\end{theorem}

We are going to further extend such results  and in particular  make good on the claim in \cite{FR21} that extensional realizability may be used to show  arithmetical conservativity of $\ac_{\ft}$ over $\izf$ (as well as $\czf$).

\begin{theorem}\label{goodman}
Let $T$ be any of the theories obtained by appending to either $\czf$ or $\izf$ some or all of the axioms $\ac_{\omega}$, $\dc$, $\dc^\dagger$, $\rdc$, $\rdc^\dagger$, $\pax$. Then:
\begin{itemize}[leftmargin=5mm]
	\item $T =_{arith}T+\Prod\Sum\I\dash\ac^\dagger+\Prod\Sum\I\dash\rdc^\dagger $;
	\item $T +\rea =_{arith}T+\rea+\Prod\Sum\W\I\dash\ac^\dagger+ \Prod\Sum\W\I\dash\rdc^\dagger$.
\end{itemize}
One can allow type $1$ parameters. Namely, we have conservativity for sentences of the form
\[  \forall f_1\in\omega^\omega\cdots \forall f_n\in\omega^\omega\, \vp(f_1,\ldots,f_n) \]
where $\vp$ is arithmetic (relative to $f_1,\ldots,f_n$). 
\end{theorem}
\begin{remark}
Of course, some of the theories to the left reduce to $T+\Prod\Sum\I\dash\ac^\dagger$ and $T+\rea+\Prod\Sum\W\I\dash\ac^\dagger$ respectively.
\end{remark}

Friedman and Scedrov's \cite{friedman_scedrov84} versions of generic realizability and forcing  only work for set theories without extensionality. We show how to adapt the argument to an extensional setting. Definition \ref{forcing} fully specifies our forcing relation. It is a routine although tedious  matter to check all the details.  By the way, one should also  prove that every theory $T$ as above  is self-validating with respect to forcing. This is not difficult to see.

It is convenient to work in a conservative extension $S$ of $T$ obtained by adding the constant symbol $\omega$ and an axiom saying that $\omega$ is the set of natural numbers.\footnote{All axioms and rules apply to the formulas of the expanded language. In the case of bounded or $\Delta_0$ separation, where the collection of formulas is syntactically defined, we  enlarge the stock of bounded formulas accordingly. These are built up from atomic formulas in the new language by closing under logical connectives and bounded quantifiers of the form $\forall x\in t$ and $\exists x\in t$, where $t$ is a term not containing $x$. Needless to say, as long as we add only constants, since parameters are allowed, bounded separation in the enriched language  is equivalent to bounded separation in the primitive language with $\in$ and $=$ only.}
Now, we can clearly identify number quantifiers with bounded quantifiers of the form $\forall x\in\omega$ and $\exists x\in\omega$.  
On the other hand, every  primitive recursive number-theoretic function and hence relation is  canonically definable by a $\Sigma$ formula in the language of set theory (with just $\in$ as primitive).\footnote{A formula is $\Sigma$ if it is built up from bounded formulas by means of $\land$, $\lor$, bounded quantifiers $\forall x\in y$ and $\exists x\in y$, and unbounded existential quantifiers $\exists x$.}
We call such formulas primitive recursive. Kleene's $T$-predicate $T(e,x,z)$ and the function $U(z)$ are among them.  We then say that a formula is arithmetic if it is built up from primitive recursive formulas  by using $\land$, $\lor$, $\neg$, $\imp$, and number quantifiers  $\forall x\in\omega$ and $\exists x\in\omega$.

\begin{proof}[Proof of Theorem \ref{goodman}]
We first consider the case without parameters. The goal is to show that for every arithmetic sentence $\vp$ of $S$, 
\begin{itemize}[leftmargin=5mm]
	\item if $S+\Phi\vdash \vp$, then $S\vdash \vp$, 
	\item if $S+\rea+\Psi\vdash \vp$, then $S+\rea\vdash \vp$,
\end{itemize}
where $\Phi=\Prod\Sum\I\dash\ac^\dagger+\Prod\Sum\I\dash\rdc^\dagger$ and $\Psi=\Prod\Sum\W\I\dash\ac^\dagger+\Prod\Sum\W\I\dash\rdc^\dagger$.\vspace{1mm}

Step 1. Let $S^g$ be the theory obtained from $S$ by adjoining a new constant symbol $g$ and an axiom saying the $g$ is a partial function from $\omega$ to $\omega$. 
Within $S^g$, one can define extensional generic realizability $a=b\fo_\omega \vp$, where $a,b\in\omega$ and $\vp$ is a formula of $T$ with parameters in  $\vo$, by using partial recursive application relative to the generic oracle $g$. To be precise, we use Kleene's relativized $T$-predicate $T(e,x,z,f)$, where $f$ is a partial function from $\omega$ to $\omega$. 
So, for example, the clause for implication reads as follows:
\begin{align*}
	a=b&\fo_\omega \vp\imp\psi && \Leftrightarrow && \forall c,d\in\omega \,  (c=d\fo_\omega \vp\imp  \{a\}^g(c)=\{b\}^g(d)\fo_\omega \psi) 
\end{align*}
For $\land$ and $\lor$, we may use a primitive recursive pairing function $(n,m)$ on natural numbers with primitive recursive projections $(n)_0$ and $(n)_1$. So, for example, we let
\begin{align*}
	a=b& \fo_\omega \vp\land \psi && \Leftrightarrow && (a)_0=(b)_0\fo_\omega \vp \land (a)_1=(b)_1\fo_\omega \psi 
\end{align*}

One can lift the realizability interpretation to every formulas $\vp$  of $S$   by interpreting $\omega$ with $\dot\omega\in\vo$. Here, 
$\dot\omega=\{\pair{n,n,\dot n}\mid n\in\omega\}$, where $\dot n=\{\pair{m,m,\dot m}\mid m\in n\}$.

One proves that 
\begin{equation} \tag{1} \text{if $S+\Phi\vdash\vp$, then   $S^g\vdash \exists e\in\omega\, (e\fo_\omega \vp)$}. 
\end{equation}

Similarly for  $S+\rea+\Psi$ and $(S+\rea)^g$.  This concludes the realizability part. \vspace{1mm}

Step 2. Next, within $S$, given any definable set $\FO$ of partial functions from $\omega$ to $\omega$ such that $0\in\FO$, one  defines a forcing interpretation $p\fof \vp$, where $p\in \FO$ and $\vp$ is a formula of $T$ with parameters in $\V(\FO)$. The universe $\V(\FO)$ is  inductively defined by the clause:
\begin{itemize}[leftmargin=5mm]
	\item if $x\subseteq\FO\times \V(\FO)$ and $\pair{q,z}\in x$ whenever  $\pair{p,z}\in x$ with $q\supseteq p$, then $x\in\V(\FO)$.
\end{itemize}

In $\V(\FO)$, we have a canonical name $\mathring g$ for the generic $g$ defined as
\[ \mathring g=\{ \pair{p,\pair{\check n,\check m}_\FO}\mid p\in \FO\land \pair{n,m}\in p\}, \]
where $\check n=\{\pair{p,\check m}\mid p\in\FO\land m\in n\}$, and $\pair{x,y}_\FO$ is some suitable pairing function in $\V(\FO)$. In particular, every natural number $n\in\omega$ has a canonical name $\check n$ in $\V(\FO)$. A canonical name $\check \omega$ for $\omega$ is  given by
$  \check\omega=\{\pair{p,\check n}\mid p\in\FO\land n\in\omega\}$. 

One can then extend the forcing interpretation $p\fof \vp$ to formulas of $S^g$  by interpreting $g$ with $\mathring g$ and $\omega$ with $\check \omega$. The basic idea however is that $p\fof g(\check n)=\check m$ iff $\pair{n,m}\in p$. Moreover, 
\begin{itemize}[leftmargin=5mm]
	\item if $\{e\}^p(a)\simeq b$, then $p\fo \{\check e\}^g(\check a)\simeq \check b$,
	\item if $p\fo  \{\check e\}^g(\check a)\simeq \check b$, then for every $q\supseteq p$, there is an $r\supseteq q$ such that $\{e\}^r(a)\simeq b$.
\end{itemize} 

Then one shows that 
\begin{equation}
	\text{if $S^g\vdash \vp$, then $S\vdash \forall p\, \exists q\supseteq p\,  (q\fof\vp)$}. \tag{2}
\end{equation}
Similarly for $S^g+\rea$ and $S+\rea$. In particular, we must have that
\[ S\vdash \forall p\, \exists q\supseteq p\, (q\fof \omega \text{ is the set of natural numbers}), \text{ and} \]
\[  S\vdash \forall p\, \exists q\supseteq p\,  (q\fof g \text{ is a partial function from $\omega$ to $\omega$}). \]

Moreover, forcing is absolute for   arithmetic sentences, namely,  
\begin{equation}
	S\vdash \forall p\, (\vp\biimp (p\fof\vp)), \tag{3} 
\end{equation}
for every arithmetic sentence $\vp$ of $S$.  This is the forcing part.

Finally, one shows that if $\vp$ is an arithmetic sentence, then there exists a forcing notion $\FO$ such that 
\begin{equation}
	S\vdash \forall p\, (p\fof \forall a,b\in\omega\, ((a=b\fo_\omega \vp)\imp \vp)). \tag{4} 
\end{equation}
\vspace{1mm}

Step 3. By putting all together, the arguments runs as follows. If $S+\Phi\vdash \vp$, with $\vp$ arithmetic, let $\FO$ be as in (4). By (1),  $S^g\vdash \exists e\in\omega\, (e\fo_\omega\vp)$.  By (2), $S\vdash \forall p\, \exists q\supseteq p\, (q\fof 
\exists e\in\omega\, (e\fo_\omega\vp))$. We now reason in $S$. By the soundness of forcing and (4), we can find $p\in \FO$ such that $p\fof \vp$.   By (3), we finally obtain $\vp$.

\begin{definition}[Forcing]\label{forcing} We define the relation $p\fof \vp$, where $p\in\FO$ and $\vp$ is a formula of $T$ with parameters in $\V(\FO)$. The atomic cases are defined by transfinite recursion. 
	\begin{align*}   
		p &\fof x\in y && \Leftrightarrow && \exists z\in\V(\FO)\, (\pair{p,z}\in y\land p\fof x=z) \\
		p & \fof x=y && \Leftrightarrow  && \forall \langle q,z\rangle \in x\, (q\supseteq p\imp  \exists r\supseteq q\, (r\fof z\in y)) \text{ and }  \\
		&&&&& \forall \langle q,z\rangle \in y\, (q\supseteq p\imp \exists r\supseteq q\, (r\fof z\in x))\\ 
		p& \fof \vp\land \psi && \Leftrightarrow && p\fof \vp \land p\fof \psi \\
		p & \fof \vp\lor\psi &&  \Leftrightarrow && p\fof \vp \lor p\fof \psi \\ 
		p &\fof \neg\vp && \Leftrightarrow && \forall q\supseteq p\,  \neg (q\fof \vp) \\
		p &\fof \vp\imp\psi && \Leftrightarrow && \forall q\supseteq p\,  (q\fof \vp\imp \exists r\supseteq q\, (r\fof \psi))   \\
		p& \fof \forall x\in y\, \vp && \Leftrightarrow && \forall \pair{q,x}\in y\, (q\supseteq p\imp \exists r\supseteq q\, (r\fof \vp)) \\
		p &\fof \exists x\in y\, \vp && \Leftrightarrow && \exists x\in\V(\FO)\, (\pair{p,x}\in y\land p\fof \vp)\\
		p& \fof \forall x\, \vp && \Leftrightarrow && \forall x\in\V(\FO)\, \forall q\supseteq p\, \exists r\supseteq q\, (r\fof \vp) \\
		p & \fof \exists x\, \vp && \Leftrightarrow && \exists x\in\V(\FO)\, (p\fof \vp) 
	\end{align*}
\end{definition} 

Forcing is monotone, that is,
if $p\fof \vp$ and $q\supseteq p$, then $q\fof \vp$.

All the ingredients are in place to carry out a  proof of (1), (2), (3), (4). Here we give some  clues. \vspace{1mm}

{\bf On absoluteness.} For every  set $x$, define by transfinite recursion the name $\check x\in\V(\FO)$ by 
$\check x=\{\pair{p,\check y}\mid p\in\FO \land y\in x\}$.
We set out  to prove that for every  arithmetic formula $\vp(x_1,\ldots,x_n)$,
\begin{equation}  
	S\vdash \forall x_1\in\omega \cdots \forall x_n\in\omega \, \forall p\, (\vp(x_1,\ldots,x_n)\biimp p\fof \vp(\check x_1,\ldots,\check x_n)). \tag{5}
\end{equation}
The absoluteness of arithmetic sentences in the sense of (3) then follows.

It is an easy exercise to show that (5) is preserved under logical connectives and number quantifiers.
It is thus sufficient to show that for every primitive recursive function $f(x_1,\ldots,x_n)$, the primitive recursive formula $\vp_f(x_1,\ldots,x_n,y)$ defining $f$ satisfies (5). This is proved by (an external) induction on the generation of $f$.

For example, suppose $f(x_1,\ldots,x_n,x)$ is defined by primitive recursion from $g(x_1,\ldots,x_n)$ and $h(x_1,\ldots,x_n,x,y)$, namely, \begin{align*}f(x_1,\ldots,x_n,0)&=g(x_1,\ldots,x_n);\\ f(x_1,\ldots,x_n,x+1)&=h(x_1,\ldots,x_n,x,f(x_1,\ldots,x_n,x)). \end{align*} 
Define $\vp_f(x_1,\ldots,x_n,x,y)$ as $\exists z\, (\psi(z) \land z(x)=y)$,  where $\psi(z)$ is 
\begin{multline*}
	\fun(z)\land \dom(z)\in\omega \land	(0\in\dom(z)\imp \vp_g(x_1,\ldots,x_n,z(0))) \\ \land   \forall n\in\dom(z)\, \forall m\in n\, (n=m+1\imp  \vp_h(x_1,\ldots,x_n,n,z(m),z(n))).
\end{multline*}
(Note  that $\vp_f$ is equivalent to a $\Sigma$ formula if  $\vp_g$ and $\vp_h$ are also $\Sigma$ definable.)  One shows that (5) holds for $\vp_f(x_1,\ldots,x_n,x,y)$, by  assuming that it holds for $\vp_g$ and $\vp_h$. The left-to-right direction of (5) requires  some standard absoluteness arguments.  
For the right-to-left direction, one also uses the fact that, for every primitive recursive function $f(x_1,\ldots,x_n)$,
\begin{align*} S&\vdash \forall x_1\in\omega\cdots \forall x_n\in\omega\,  (\exists y\in\omega \, \vp_f(x_1,\ldots,x_n,y) \\ 
	& \qquad \land \forall y_0\, \forall y_1\, (\vp_f(x_1,\ldots,x_n,y_0)\land \vp_f(x_1,\ldots,x_n,y_1)\imp y_0=y_1)), 
\end{align*}
and thereby this is forced in $S$.\vspace{1mm}

{\bf On forcing self-realizability.} To obtain (4), one shows that for every arithmetic formula $\vp(x_1,\ldots,x_n)$  there is a $\FO$ such that $S$ proves: 
\begin{enumerate}[label=(\roman*), leftmargin=*]
	\item there is  $e\in\omega$ such that for every $p$ there is a $q\supseteq p$ forcing
	\begin{align*}
		& \forall a_1\in\omega\cdots\forall a_n\in\omega\, (\vp(a_1,\ldots,a_n)\imp \\ & \qquad\qquad\qquad \exists b\in\omega\, (\{\check e\}^g(a_1,\ldots,a_n,0)=b\land b\fo_\omega \vp(\dot a_1,\ldots,\dot a_n))); 
	\end{align*}
	\item for every $p$ there is a $q\supseteq p$ forcing
	\[ \forall a_1\in\omega\cdots\forall a_n\in\omega\, \forall a,b\in\omega\, ((a=b\fo_\omega \vp(\dot a_1,\ldots,\dot a_n))\imp \vp(a_1,\ldots,a_n)). \]	
\end{enumerate}
Here, $\{e\}^f(a,b)$ stands for $\{\{e\}^f(a)\}^f(b)$, and so on and so forth. On a technical note,  the dummy argument for $0$ is used to handle the case where $\vp$ is closed. 

Now, let $\vp(x_1,\ldots,x_n)$ be arithmetic. We define $\FO$ as follows. First, let $\mathbb{Q}$ be the set of partial functions $p$ from $\omega$ to $\omega$ such that $\dom(p)$ is an image of a natural number. The point of this definition is just to make sure that, within $\czf$, 
the class $\mathbb{Q}$ exists as a set and the domain of every $p\in\mathbb{Q}$ is decidable, in the sense that 
\[  \forall p\in\mathbb{Q}\, \forall x\in\omega\, (x\in\dom(p)\lor x\notin\dom(p)). \]
Note that if we are working in $\izf$, we can directly define the set $\mathbb{Q}$ of all partial functions from  $\omega$ to $\omega$ with decidable domain. 
Enumerate all subformulas $\vp_1,\ldots,\vp_k$ of $\vp$ starting from the primitive recursive ones. Let $\FO$ be the set of $p\in\mathbb{Q}$ such that for every subformula $\vp_i(y_1,\ldots, y_{j}) $ of the form $\psi_0\lor\psi_1$, for every $a_1,\ldots,a_j\in\omega$, if $a=(i,a_1,\ldots,a_j)\in\dom(p)$, then $p(a)=0\land \psi_0(a_1,\ldots,a_j)$ or $p(a)=1\land\psi_0(a_1,\ldots,a_j)$, and for every subformula $\vp_i(y_1,\ldots,y_j)$ of the form $\exists x\, \psi(x,y_1,\ldots,y_j)$, for every $a_1,\ldots,a_j\in\omega$, if $a=(i,a_1,\ldots,a_j)\in\dom(p)$, then $\psi(p(a),a_1,\ldots,a_j)$.

One can write down a formula defining $\FO$ as there are standard finitely many subformulas of $\vp$. Here, $(a_1,\ldots,a_j)$ is a primitive recursive coding of sequences of natural numbers. One then shows (i) and (ii) for every subformula of $\vp$. 

If $\vp$ is primitive recursive, then one can find a number $e$ such that
\begin{multline*} S^g\vdash  \forall a_1,\ldots,a_n\in\omega (\vp(a_1,\ldots,a_n)  \biimp \{e\}^g(a_1,\ldots,a_n)\fo_\omega \vp(\dot a_1,\ldots,\dot a_n)), 
\end{multline*}
where we identify $e$ with the corresponding numeral. This fact can be gleaned from \cite[Chapter 4, Theorem 2.6]{M84} in the case of $\izf$. The same holds for $\czf$ \cite[Proposition 8.5]{R06}.  We can then use the soundness of forcing  to get (i) and (ii). 

The index of complex subformulas is computed according to the following table. Let $i$ be the number of the given subformula according to our enumeration. Write $\vec{a}$ for the tuple $a_1,\ldots,a_j$.

\begin{align*}
	\{e_{\vp\land\psi}\}^f(\vec a,0)&\simeq (\{e_\vp\}^f(\vec a,0),\{e_\psi\}^f(\vec a,0))   \\
	\{e_{\vp\lor\psi}\}^f(\vec a,0)&\simeq \begin{cases}
		\{e_\vp\}^f(\vec a,0) & \text{if } f(i,\vec a)=0 \\
		\{e_{\psi}\}^f(\vec a,0)& \text{if } f(i,\vec a)=1 \end{cases} \\
	\{e_{\forall x\in\omega\, \vp(x)}\}^f(\vec a,0,b)&\simeq \{e_\vp\}^f(\vec a,b,0)\\
	\{e_{\exists x\in\omega\, \vp(x)}\}^f(\vec a,0)&\simeq (f(i,\vec a),\{e_\vp\}^f(\vec a,f(i,\vec a),0))
\end{align*}

By way of example, let us illustrate how to obtain (i)  in the existential case. For simplicity, suppose we have a subformula of the form $\exists y\in\omega \, \vp(x,y)$. We work in $S$. We have already found an index $e_\vp$ for $\vp$. So, let $e$ be such that 
\[ \{e\}^f(a,0)\simeq (f(i,a),\{e_\vp\}^f(a,f(i,a),0)). \] 

Let $p\in \FO$. We search for a $q\supseteq p$ such that
\[  q\fof \forall a\in\omega\, (\exists y\in\omega\, \vp(a,y)\imp  \{\check e\}^g(a,0)\fo_\omega \exists y\in\omega\, \vp(\dot a,y)). \] 
Notice that, within $S^g$, 
\[ (a=b\fo_\omega \exists x\in\omega\, \vp(x)) \biimp \exists n\in\omega\, ((a)_0=(b)_0=n\land 
(a)_1=(b)_1\fo_\omega \vp(\dot n)), \]
and thereby this is forced in $S$. It is then enough to find an extension $q\supseteq p$ such that
\begin{align*} 
	q&\fof  \forall a\in\omega\, (\exists y\in\omega\, \vp(a,y)\imp \exists c\in\omega\, \exists b\in\omega\, (\{\check e\}^g(a,0)\simeq c  \\ \tag{1} &\qquad\qquad\qquad\qquad\qquad\qquad\qquad \land (c)_0=b  \land (c)_1\fo_\omega \vp(\dot a,\dot b))). 
\end{align*}
We claim that $p$ already satisfies (1). Let $a\in\omega$, $q\supseteq p$,  and suppose 
$q\fof \exists y\in\omega\, \vp(\check a,y)$.  We are clear if we find $r\supseteq q$ such that
\begin{align*} 
	r&\fof \exists c\in\omega\, \exists b\in\omega\, (\{\check e\}^g(\check a,0)\simeq c\land (c)_0=b \\ & \tag{2} \qquad\qquad \land \forall x\, (i(\check a,x)\imp (c)_1\fo_\omega \vp(x,\dot b))),
\end{align*}
where the formula $i(x,y)$ expresses that $y=\dot x$ for $x\in\omega$. By definition, there is a $b\in\omega$ such that $q\fof \vp(\check a,\check b)$. By absoluteness, $\vp(a,b)$ holds. Now, we can decide whether $(i,a)\in\dom(q)$. If not, we can take $q'=q\cup\{\pair{(i,a),b}\}$. By the  assumption on $e_\vp$, there is an $s\supseteq q'$ such that  \begin{align*} 
	& s\fof  \exists c_1\in\omega\, (\{\check e_\vp\}^g(\check a,\check b,0)\simeq c_1 \\
	& \tag{3} \qquad\qquad\qquad \land \forall x\, \forall y\, (i(\check a,x)\land i(\check b,y) \imp  c_1\fo_\omega \vp(x,y))).
\end{align*}
One can see that 
$s\fof \{\check e\}^g(\check a,0)\simeq   (\check b,\{\check e_\vp\}^g(\check a,\check b,0))$.  
By  (3) and soundness, one can finally find an $r\supseteq s$ satisfying (2), as desired.\vspace{1mm}

{\bf Adding parameters.} We briefly discuss how to generalize such conservation result to arithmetic formulas with type $1$ parameters. For example, suppose $S+\Phi\vdash \forall f\in\omega^\omega\, \vp(f)$. Then $(S+\Phi)_f\vdash \vp(f)$,
where, in general, $T_f$ is obtained from $T$ by adjoining a constant symbol $f$ and an axiom saying that $f$ is a function from $\omega$ to $\omega$.  It clearly suffices to show that $S_f\vdash \vp(f)$. If we do so, then   $S\vdash \forall f\in\omega^\omega\, \vp(f)$

Now consider the theory $(S_f)^g$ obtained  from $S_f$ by further extending the language with a constant symbol $g$ and an axiom saying that $g$ is a partial function from $\omega$ to $\omega$. Within $(S_f)^g$, consider extensional generic realizability $a=b\fo_\omega \vp$, where $a,b\in\omega$ and $\vp$ is a formula of $S_f$, by using partial recursive application relative to the oracle $f\oplus g=\{\pair{2n,f(n)}\mid n\in\omega\}\cup\{\pair{2n+1,g(n)}\mid n\in\dom(g)\}$. Use the name $\dot f=\{\pair{n,n,\pair{\dot n,\dot m}_\omega}\mid f(n)=m\}\in\vo$ to define realizability on formulas containing $f$. Note that we could extend the interpretation to every formula of $(S_f)^g$ by interpreting $g$ with $\dot g$, which is defined in a similar way as $\dot f$.
Show that  
\[  \text{if } (S+\Phi)_f\vdash \psi, \text{ then }   (S_f)^g\vdash \exists e\in\omega\, (e\fo_\omega \psi). \] 
Similarly for $(S+\rea+\Psi)_f$ and   $(S_f+\rea)^g$.  In particular, we must have
\[  (S_f)^g\vdash \exists e\in\omega\, (e\fo_\omega f\colon\omega\to\omega). \] 
This is the realizability part. 

In $S_f$, given any definable set $\FO$ of partial functions from $\omega$ to $\omega$, where now we can use the parameter $f$ in defining $\FO$, we extend the forcing interpretation $p\fof \vp$ of Definition \ref{forcing} to formulas of $(S_f)^g$ by interpreting  $g$ with $\mathring g$, like before, and $f$ with $\check f$. Note that
$\check f=\{\pair{p,\check z}\mid p\in\FO\land z\in f\}=\{\pair{p,\pair{\check n,\check m}_\FO}\mid p\in\FO\land f(n)=m\}$. 
The rest of the argument is  roughly as before, but with everything relativized to the parameter $f$. In particular, one must show that arithmetic (in $f$) formulas are absolute for forcing. Finally, given $\vp(f)$, one defines $\FO$ as before, but relative to $f$, and show that in $\V(\FO)$, $\vp(f)$ is self-realizing. 
\end{proof}

\appendix
\section{}\label{app}
The proofs of Theorem \ref{realizers} and Theorem \ref{canonical} will rely on the following two facts.

\begin{theorem}[$\czf$]\label{def}
There is a formula $\vartheta(u,Y)$, such that 
\[  \Prod\Sum\I=\{Y\mid \exists u\, \vartheta(u,Y)\}\ \ \text{and}\ \  \vartheta(u,Y)\biimp \bigvee_{i=0}^4 \vartheta_{i}(u,Y), \]
where 
\begin{align*}
	& \vartheta_0(u,Y) && \Leftrightarrow && Y\in\omega,  \\
	&\vartheta_1(u,Y) && \Leftrightarrow && Y=\omega,  \\
	&\vartheta_2(u,Y) && \Leftrightarrow &&  \exists X\, \exists F\, (\psi(u,X,F)\land Y=\Prod(X,F)),      \\
	&\vartheta_3(u,Y) && \Leftrightarrow && \exists X\, \exists F\, (\psi(u,X,F)\land Y=\Sum(X,F)), \\
	& \vartheta_4(u,Y) &&  \Leftrightarrow && \exists X\, (\exists v\in u\, \vartheta(v,X)\land \exists x\in X\, \exists y\in X\, (Y=\I(x,y))), 
\end{align*}
with $\psi(u,X,F)$ being a shorthand for
\[ \exists v\in u\, \vartheta(v,X)\land  \fun(F) \land \dom(F)=X\land   \forall x\in X\, \exists v\in u\, \vartheta(v,F(x)). \] 
Similarly, there is a formula  $\vartheta^*(u,Y)$, such that 
\[  \Prod\Sum\W\I=\{Y\mid \exists u\, \vartheta^*(u,Y)\} \ \ \text{and}\ \ 
\vartheta^*(u,Y)\biimp \bigvee_{i=0}^5 \vartheta^*_{i}(u,Y), \]
where  $\vartheta^*_i$ for $i=0,\ldots,4$ is defined in a similar fashion as $\vartheta_i$, and 
\begin{align*}
	&\vartheta^*_5(u,Y) && \Leftrightarrow && \exists X\, \exists F\, (\exists v\in u\, \vartheta^*(v,X)\land \fun(F) \land \dom(F)=X\land \\
	&&&&&   \forall x\in X\, \exists v\in u\, \vartheta^*(v,F(x))\land Y=\W(X,F)).
\end{align*}	
\end{theorem}

\begin{theorem}[$\czf+\rea$]\label{defW}
If $F$ is a function with domain $X$, then for every set $Y$, 
\begin{align*}
	Y&=\W(X,F) \biimp \forall y\, (y\in Y\biimp \exists x\in X\, \exists f\, (y=\pair{x,f}\land \fun(f)\land \\
	& \qquad \qquad \dom(f)=F(x)\land \forall u\in\dom(f)\, (f(u)\in Y))).
\end{align*}
\end{theorem}

\begin{proof}[Proof of Theorem \ref{realizers}]
From now on, let
\[ F=F_{\sigma,i}=\{\pair{a,b,\vpair{a^\sigma,X_\tau}}\mid a\sim_\sigma b\land ia=\tau\}. \] 

(i) A closed term $\mb o$ such that $\mb o\fo X_{\sf N}=\omega$ is given in the proof of Theorem \ref{sound}. In fact, the set $X_{\sf N}=\dot \omega$ is the canonical name for $\omega$ used to realize infinity. \vspace{1mm}

(ii) We aim for 
$\mb{fun}\fo \fun(F)\land \dom(F) = X_\sigma$. 
It suffices to look for closed terms $\mb r,\mb f, \mb t$ such that: 
\begin{align*}  \mb r&\fo \forall z\in F\, \exists x\in X_\sigma\, \exists y\, \op(z,x,y)  \quad (\text{namely, } F\text{ is a relation }\land  \dom(F)\subseteq X_\sigma),  \\
	\mb t& \fo \forall x\in X_\sigma\, \exists y\, \exists z\in F\, \op(z,x,y) \quad (\text{namely, } X_\sigma\subseteq \dom(F)), \\
	\mb f &\fo F \text{ is functional.}  	& 
\end{align*}
We can easily arrange for $\mb r$ and $\mb t$. Let us find $\mb f$.  By unraveling the definition,  we want $\mb f$ such that, if $a\sim_\sigma b $, $\breve a\sim_\sigma \breve b$ and
\[  c=d\fo \op(\vpair{a^\sigma,X_{\tau}},x,y_0)\land \op(\vpair{\breve a^\sigma,X_{\breve\tau}},x,y_1), \]
where $ia=\tau$ and $i\breve a=\breve \tau$, then $\mb fa\breve ac=\mb fb\breve bd\fo y_0=y_1$. 
By the properties of pairing, from $c$ and $d$ as above we get  $\fo a^\sigma=\breve a^\sigma$. By Lemma \ref{inj}, we must have $a\sim_\sigma \breve a$. Since $i$ is a family of types over $\sigma$, $ia\sim i\breve a$. That is, $\tau\sim \breve\tau$. By Lemma \ref{inv}, $X_\tau=X_{\breve\tau}$.  As in the proof of Theorem \ref{canonical choice}, by the properties of equality, there is a closed term $\mb q$ such that
$\mb q c=\mb q d\fo y_0=y_1$.
We  may let  $\mb f=_{\mathrm{def}}\lambda a\breve ac. \mb q c$.\vspace{1mm}

(iii)  We aim for
$\mb{prod}\fo  X_{\Prod_\sigma i}=\Prod(X_\sigma, F)$. 
We abbreviate 
\begin{align*}
	&\vartheta_r(f,X,F) && \Leftrightarrow &&  f\subseteq \Sum(X,F)\\
	&&&\Leftrightarrow &&   \forall z\in f\, \exists x\in X\, \exists y\, \exists Z\in F\, \exists Y\, (\op(z,x,y)\land \op(Z,x,Y)\land y\in Y), \\
	&\vartheta_t(X,f)  && \Leftrightarrow &&  X\subseteq\dom(f) \\
	&&&\Leftrightarrow &&  \forall x\in X\, \exists z\in f\, \exists y\, \op(z,x,y), \\
	&\vartheta_f(f)  && \Leftrightarrow && f \text{ is functional}\\
	&&& \Leftrightarrow && \forall z_0\in f\, \forall z_1\in f\, \forall x\, \forall y_0\, \forall y_1\, (\op(z_0,x,y_0)\land \op(z_1,x,y_1)\imp y_0=y_1).
\end{align*}
Let $\alpha=\Prod_\sigma i$. For the direction left to right, it suffices to look for closed terms $\mb r$, $\mb t$, $\mb f$ such that whenever $f\sim_\alpha g$, then 
\begin{align*} \tag{1} \mb rf=\mb rg &\fo \vartheta_r(f^\alpha,X_\sigma,F), \\
	\tag{2} \mb tf=\mb tg &\fo \vartheta_t(X_\sigma,f^\alpha),\\
	\tag{3}  \mb ff=\mb fg &\fo \vartheta_f(f^\alpha). 
\end{align*}
For the direction right to left, we look for a closed term $\mb c$ such that, for every $h\in\va$, if
$c=d\fo \vartheta_r(h,X_\sigma,F)\land \vartheta_t(X_\sigma,h)\land \vartheta_f(h)$,
then 
$\mb c c=\mb c d\fo h\in X_{\Prod_\sigma i}$.

Let us start with the forward direction. 

(1) Let $f\sim_\alpha g$,  $a\sim_\sigma b$, $ia=\tau$ and $fa=e$. We want
\[   \mb rfa=\mb rgb\fo  \exists x\in X_\sigma\, \exists y\, \exists Z\in F\, \exists Y\, (\op(\vpair{a^\sigma,e^\tau},x,y)\land \op(Z,x,Y)\land y\in Y).\]
We let $(\mb rfa)_0=a$ and set out for
\[  (\mb rfa)_1=(\mb rgb)_1\fo   \exists Z\in F\, \exists Y\, (\op(\vpair{a^\sigma,e^\tau},a^\sigma,e^\tau)\land \op(Z,a^\sigma,Y)\land e^\tau\in Y).\]
We let $(\mb rfa)_{10}=a$, $Y=X_\tau$ and set out for 
\[  (\mb rfa)_{11}=(\mb rgb)_{11}\fo \op(\vpair{a^\sigma,e^\tau},a^\sigma,e^\tau)\land \op(\vpair{a^\sigma,X_\tau},a^\sigma,X_\tau)\land e^\tau\in X_\tau.\]
Let $(\mb rfa)_{110}=(\mb rfa)_{1110}=\mb v$, where $\mb v\fo \op(\vpair{x,y},x,y)$ exists by  Lemma \ref{pairs}. Let $(\mb rfa)_{1111}=\mb p(ca)\mb i$, where $\mb i\fo \forall x\, (x=x)$.  Then 
$\pair{fa,gb,e^\tau}\in X_\tau$ and 
\[  \mb r=_{\mathrm{def}}\lambda f\lambda a. \mb p a(\mb p a(\mb p \mb v(\mb p \mb v(\mb (fa)\mb i))))) \]
is as desired. 

(2)  Exercise.

(3) Let $f\sim_\alpha g$. Let  $a_j\sim_\sigma b_j$, $ia_j=\tau_j$, $fa_j=e_j$ for $j=0,1$.  We wish to find $\mb f$ such that, if $c=d\fo  \op(\vpair{a_0^\sigma,e_0^{\tau_0}},x,y_0)\land \op(\vpair{a_1^\sigma,e_1^{\tau_1}},x,y_1)$, then 
$\mb ffa_0a_1c=\mb fgb_0b_1d\fo y_0=y_1$, for  all $x,y_0,y_1\in\va$. Suppose we have $c$ and $d$ as above.  Then in particular $\fo a_0^\sigma=a_1^\sigma$.
By Lemma \ref{inj}, it follows that $a_0\sim_\sigma a_1$, and hence $\tau_0\sim\tau_1$ and $e_0\sim_{\tau_0} e_1$. Thus by Lemma \ref{inj} again, $y=e_0^{\tau_0}=e_1^{\tau_1}$. Then $c=d\fo  \op(\vpair{a_0^\sigma,y},x,y_0)\land \op(\vpair{a_1^\sigma,y},x,y_1)$.  By soundness and the properties of pairing, $\mb q c=\mb qd\fo y_0=y_1$, for some closed term $\mb q$.  
Thus 
$\mb f=_{\mathrm{def}}\lambda f\lambda a_0\lambda a_1\lambda c. \mb qc$
is as desired.

Conversely, let $h\in\va$ and suppose 
\begin{align*}   c_0=d_0 &\fo \vartheta_r(h,X_\sigma,F) & (\text{namely, } h\subseteq\Sum(X_\sigma,F)),    \tag{4} \\
	c_{10}=d_{10} &\fo \vartheta_t(X_\sigma,h) & (\text{namely, } X_\sigma\subseteq\dom(h)),
	\tag{5}\\
	c_{11}=d_{11} &\fo \vartheta_f(h) & (\text{namely, } h\text{ is functional}).	\tag{6}
\end{align*}
We wish to find $\mb c$ such that   
$\mb c c=\mb c d\fo h\in X_{\Prod_\sigma i}$.
We thus aim for $f=(\mb cc)_0\sim_\alpha (\mb cd)_0=g$ and $(\mb cc)_1=(\mb cd)_1\fo h=f^\alpha$.
Let us first find $f$ and $g$. If $a\sim_\sigma b$, then  by (5) there are $z,y\in\va$ such that
\begin{equation}   
	\pair{(c_{10}a)_0,(d_{10}b)_0,z}\in h,  \tag{7}
\end{equation}
\begin{equation}  (c_{10}a)_1=(d_{10}b)_1\fo \op(z,a^\sigma,y). \tag{8} 
\end{equation}
From (4) and (7), it follows  that
\[   \breve{a}=(c_0(c_{10}a)_0)_0\sim_\sigma (d_0(d_{10}b)_0)_0=\breve{b}, \]
\[   \tilde{a}=(c_0(c_{10}a)_0)_{10}\sim_\sigma (d_0(d_{10}b)_0)_{10}=\tilde{b}, \]
and for some $\breve{y},Y\in\va$,
\begin{equation}   (c_0(c_{10}a)_0)_{11}= (d_0(d_{10}b)_0)_{11}\fo \op(z,\breve{a}^\sigma,\breve{y})\land \op(\vpair{\tilde{a}^\sigma,X_{\tilde{\tau}}},\breve{a}^\sigma,Y)\land \breve{y}\in Y, \tag{9} 
\end{equation}
where  $\tilde\tau=i \tilde{a}$. Now, by (8) and (9) we obtain
$\fo a^\sigma=\breve{a}^\sigma$ and  $\fo \tilde{a}^\sigma=\breve{a}^\sigma$,
and hence by Lemma \ref{inj} we have $a\sim_\sigma \breve{a}\sim_\sigma \tilde{a}$. Therefore,   $a^\sigma=\breve a^\sigma$, $\tau=ia\sim i\breve{a}=\breve{\tau}$,  and $X_\tau=X_{\breve{\tau}}$.  By the properties of  pairing and equality, from (8) and (9) we obtain
$\mb q ca=\mb q db\fo y\in X_\tau$, 
for some closed term $\mb q$. Hence
\begin{equation}   
	e=(\mb q ca)_0\sim_\tau (\mb qdb)_0, \tag{10} 
\end{equation}
\begin{equation} (\mb qca)_1=(\mb qdb)_1\fo y=e^\tau.  \tag{11}  
\end{equation}

Therefore, by (10), by letting $(\mb c c)_0=\lambda a.(\mb qca)_0$, we have 
\[  f=(\mb cc)_0\sim_\alpha (\mb cd)_0=g. \]
It remains to show 
$(\mb c c)_1=(\mb c d)_1\fo  h=f^\alpha$.

For the right to left inclusion, let $a\sim_\sigma b$, $e=fa$ and $\tau=ia$. We want 
\[  ((\mb c c)_1a)_1= ((\mb c d)_1b)_1 \fo \vpair{a^\sigma,e^\tau}\in h. \]
By the foregoing discussion, we have that (7), (8) and (11) hold true for some  $z,y\in\va$.
By the properties of  pairing and equality, it follows from (8) and (11) that  
$\mb rca=\mb rdb\fo \vpair{a^\sigma,e^\tau}=z$,  
for some closed term $\mb r$. We may then let  
$((\mb c c)_1a)_{1}=\mb p (c_{10}a)_0 (\mb rca)$.

The left to right inclusion is treated in a similar manner.  Suppose $\pair{\breve c,\breve d,\breve z}\in h$. By (4), 
$a=(c_0\breve c)_0\sim_\sigma (d_0\breve d)_0=b$, 
and for some $\breve y\in\va$,
$(c_0\breve c)_{110}= (d_0\breve d)_{110}\fo \op(\breve z,a^\sigma,\breve y)$. 
As before, there are $z,y\in\va$ such that  (7), (8), and (11) hold true of $a$ and $b$.   On account of (6),  it follows from (7), (8), (11) that
$\mb l c\breve c=\mb ld\breve d\fo \breve z=\vpair{a^\sigma,(fa)^\tau}$, 
for some closed term $\mb l$.
Then $\mb c$ such that
$((\mb c c)_1\breve c)_0=\mb p (c_0\breve c)_0 (\mb lc\breve c)$
does the job.\vspace{1mm}

(iv)   We aim for 
$\mb{sum}\fo X_{\Sum_\sigma i}=\Sum(X_\sigma, F)$.  
Let $\beta=\Sum_\sigma i$. For the direction from left to right, we ask for a closed term $\mb l$ such that 
\[  \mb l\fo \forall z\in X_{\Sum_\sigma i}\, \exists x\in X_\sigma\, \exists y\, \exists Y\, \exists Z\in F\, (\op(z,x,y)\land \op(Z,x,Y)\land y\in Y). \]
Let $a\sim_\beta b$. We want
\[ \mb la=\mb l b\fo \exists x\in X_\sigma\, \exists y\, \exists Y\, \exists Z\in F\, (\op(\vpair{a_0^\sigma,a_1^\tau},x,y)\land \op(Z,x,Y)\land y\in Y), \]
where $\tau=ia_0$. By definition, $a_0\sim_\sigma b_0$. Let $(\mb l a)_0=a_0$. Then we aim for
\[ (\mb la)_1=(\mb lb)_1\fo \exists Z\in F\, (\op(\vpair{a_0^\sigma,a_1^\tau},a_0^\sigma,y)\land \op(Z,a_0^\sigma,Y)\land y\in Y), \]
for some $y,Y\in\va$. Let $y=a_1^\tau$, and $Y=X_\tau$. Note that 
$\pair{a_0,b_0,\vpair{a_0^\sigma,X_\tau}}\in F$.  
Set $(\mb la)_{10}=a_0$. We thus want
\[ (\mb la)_{11}=(\mb lb)_{11}\fo \op(\vpair{a_0^\sigma,a_1^\tau},a_0^\sigma,a_1^\tau)\land \op(\vpair{a_0^\sigma,X_\tau},a_0^\sigma,X_\tau)\land a_1^\tau\in X_\tau. \]
By the properties of pairing, we just need to worry about the last conjunct. Now, $a_1\sim_\tau b_1$, and thereby
$\pair{a_1,b_1,a_1^\tau}\in X_\tau$.
Therefore,
$\mb p a_1\mb i=\mb p b_1\mb i\fo  a_1^\tau\in X_\tau$,
where $\mb i\fo \forall x\, (x=x)$.

We now consider the converse direction. We look for a closed term $\mb r$ such that
\[ \mb r\fo \forall z\, \forall x\in X_\sigma\, \forall y\, \forall Z\in F\, \forall Y\, (\op(z,x,y)\land \op(Z,x,Y)\land y\in Y\imp z\in X_{\Sum_\sigma i}). \]
Let $z,y\in\va$ and $\breve a\sim_\sigma\breve b$. We want
\[ \mb r\breve a=\mb r\breve b\fo \forall Z\in F\, \forall Y\, (\op(z,\breve a^\sigma,y)\land \op(Z,\breve a^\sigma,Y)\land y\in Y\imp z\in X_{\Sum_\sigma i}). \] 
Let $a\sim_\sigma b$ and $\tau=ia$. We set out for 
\[ \mb r\breve aa=\mb r\breve bb\fo \op(z,\breve a^\sigma,y)\land \op(\vpair{a^\sigma,X_\tau},\breve a^\sigma,Y)\land y\in Y\imp z\in X_{\Sum_\sigma i}, \]
for every $Y\in\va$. Fix a $Y\in\va$ and suppose 
\[ c=d\fo \op(z,\breve a^\sigma,y)\land \op(\vpair{a^\sigma,X_\tau},\breve a^\sigma,Y)\land y\in Y. \]
By the properties of pairing and equality, we have $\fo a^\sigma=\breve a^\sigma$ and
\[  \mb q c=\mb qd\fo \op(z,a^\sigma,y)\land y\in X_\tau, \]
for some closed term $\mb q$.  By Lemma \ref{inj}, we then obtain $a^\sigma=\breve a^\sigma$. Therefore, 
\[    \mathring c=(\mb q c)_{10}\sim_\tau (\mb q d)_{10}=\mathring d, \]
\[    (\mb q c)_{11}\sim_\tau (\mb q d)_{11} \fo y= \mathring c^\tau. \]
By the usual properties of pairing and equality, we have
$\mb z c=\mb z d\fo z=\vpair{a^\sigma,\mathring c^\tau}$,
for some closed term $\mb z$. Now, 
$\mb p a\mathring c\sim_\beta\mb p b\mathring d$,
and so 
$\pair{\mb p a\mathring c,\mb p b\mathring d, \vpair{a^\sigma,\mathring c^\tau}}\in X_{\Sum_\sigma i}$.
Therefore, 
\[ \mb p(\mb p a(\mb q c)_{10})(\mb z c)=\mb p(\mb p b(\mb q d)_{10})(\mb zd)\fo z\in X_{\Sum_\sigma i}. \]
Then $\mb r$ such that $\mb r\breve a ac=\mb p(\mb p a(\mb q c)_{10})(\mb z c)$ does the job.\vspace{1mm}

(v) We aim for
$\mb{id}\fo X_{\I_\sigma(a,b)}=\I(a^\sigma,b^\sigma)$. 
For the left to right direction, we want a closed term $\mb l$ such that whenever $\pair{c,d,z}\in X_{\I_\sigma(a,b)}$, then
$\mb{l}c=\mb ld\fo \nu(z)\land a^\sigma=b^\sigma$, where $\nu(z)$ expresses that $z=0$. 
Suppose $\pair{c,d,z}\in X_{\I_\sigma(a,b)}$. By definition, $c=d=\mb 0$, $z=0$, and $a\sim_\sigma b$. It is easy to find $\mb t$ such that $\mb t\fo \nu(0)$. On the other hand, by Lemma \ref{inj}, from $a\sim_\sigma b$ we obtain $a^\sigma=b^\sigma$. Let 
$\mb l=_{\mathrm{def}}\lambda c.\mb p\mb q\mb i$, 
where $\mb i\fo \forall x\, (x=x)$.

For the right to left direction, we aim for a closed term $\mb r$ such that for every $z\in\va$,
$\mb r\fo \nu(z)\land a^\sigma=b^\sigma \imp z\in X_{\I_\sigma(a,b)}$.
Note that  $\fo \nu(z)$ implies $z=0$. By Lemma \ref{inj}, $\fo a^\sigma=b^\sigma$ implies $a\sim_\sigma b$. Thus let $\mb r=_{\mathrm{def}}\lambda c.\mb p\mb 0\mb i$,
where $\mb i$ is as before. \vspace{1mm}

(vi) We aim for 
$\mb w\fo X_{\W_\sigma i}=\W(X_\sigma,F)$.
Set $\delta=\W_\sigma i$.
By Theorem \ref{defW}, it suffices to find closed terms $\mb l$ and $\mb r$ such that 
\begin{align*} \mb l &\fo \forall y\in X_\delta\, \exists x\in X_\sigma\, \exists f\, \vartheta(y,x,f), \\
	\mb r &\fo \forall y\, (\exists x\in X_\sigma\, \exists f\, \vartheta(y,x,f)\imp y\in X_\delta),
\end{align*}
where
$\vartheta(y,x,f) \Leftrightarrow   \op(y,x,f)\land \vartheta_0(x,f)\land \vartheta_1(f)\land \vartheta_2(f)$, 
with
\begin{align*}
	&\vartheta_0(x,f) && \Leftrightarrow &&  f \text{ is a relation}\land \dom(f)=F(x) \\
	& && \Leftrightarrow &&
	\forall Z\in F\, \forall X\, (\op(Z,x,X)\imp \forall v\in f\, \exists u\in X\, \exists z\, \op(v,u,z)\land  \\
	&&&&& \qquad\qquad\qquad\qquad\qquad\qquad\qquad \forall u\in X\, \exists z\, \exists v\in f\, \op(v,u,z)),\\
	&\vartheta_1(f) && \Leftrightarrow && f \text{ is functional} \\
	&&& \Leftrightarrow &&  \forall v_0\in f\, \forall v_1\in f\, \forall u\, \forall z_0\, \forall z_1\, (\op(v_0,u,z_0)\land \op(v_1,u,z_1)\imp z_0=z_1), \\
	&\vartheta_2(f) && \Leftrightarrow &&   \ran(f)\subseteq X_\delta      \\
	&&& \Leftrightarrow &&  \forall v\in f\, \forall u\, \forall z\, (\op(v,u,z)\imp z\in X_\delta).
\end{align*}
Remember that 
\begin{align*}
	F&=F_{\sigma,i}=\{\pair{a,b,\vpair{a^\sigma,X_{ia}}}\mid a\sim_\sigma b\}, \\
	X_\delta&=\{\pair{c,d,c^\delta}\mid c\sim_\delta d\}, \\
	c^\delta&=\vpair{c_0^\sigma,f_c}, \ \ \text{where }
	f_c=\{\pair{p,q,\vpair{p^{ic_0},(c_1p)^\delta}}\mid p\sim_{ic_0} q\}, \\
	c\sim_{\delta} d&\biimp c_0\sim_\sigma d_0\land \forall p\, \forall q\, (p\sim_{ic_0}q\imp c_1p\sim_\delta d_1q).
\end{align*}

From left to right, note that if $c\sim_\delta d$, then $c_0\sim_\sigma d_0$. Let $(\mb lc)_0=c_0$. We then want
$(\mb lc)_1=(\mb ld)_1\fo \exists f\, \vartheta(c^\delta,c_0^\sigma,f)$. 
The obvious candidate for $f$ is $f_c$. The details are tedious but straightforward.

Let us consider the slightly more interesting direction. Let  $y\in\va$ and  suppose
\begin{equation} 
	a=b\fo \exists x\in X_\sigma\, \exists f\, \vartheta(y,x,f). \tag{16} 
\end{equation}
The goal is to find $\mb r$ such that $\mb ra=\mb rb\fo y\in X_\delta$.
By undoing (16), we obtain 
\begin{equation}  
	a_0\sim_\sigma b_0\ \ \text{and}\ \  a_1=b_1\fo \vartheta(y,a_0^\sigma,f), \tag{17} 
\end{equation}
for some $f\in\va$. It follows from the second half of (17) that
$a_{10}=b_{10}\fo \op(y,a_0^\sigma,f)$. Also, we can find  closed terms $\mb q_0,\mb q_1,\mb q_2$ such that  $\mb q_0a=\mb q_0b\fo \vartheta_0(a_0^\sigma,f)$ and  $\mb q_ia=\mb q_ib\fo \vartheta_i(f)$ for $i=1,2$.
Because $a_0\sim_\sigma b_0$ by the first half of (17),  we can apply $\mb q_0$ and  obtain 
\begin{equation}    \mb q_0aa_0=\mb q_0bb_0\fo \forall X\, (\op(\vpair{a_0^\sigma,X_{ia_0}},a_0^\sigma,X)\imp \psi_0(f,X)), \tag{18} 
\end{equation}
where
\begin{align*}
	&\psi_0(f,X) && \Leftrightarrow && \forall v\in f\, \exists u\in X\, \exists z\, \op(v,u,z)\land  \forall u\in X\, \exists z\, \exists v\in f\, \op(v,u,z).
\end{align*}
By the properties of pairing with $X=X_{ia_0}$, it follows from (18) that
\begin{equation}  \mb ta=\mb t b\fo \psi_0(f,X_{ia_0}), \tag{19} 
\end{equation}
for some closed term $\mb t$. Suppose $p\sim_{ia_0} q$. By (19),
$(\mb ta)_1p=(\mb t b)_1q\fo \exists z\, \exists v\in f\, \op(v,p^{ia_0},z)$,
and hence for some  $z,v\in\va$,
\begin{equation}   \pair{((\mb ta)_1p)_0, ((\mb t b)_1q)_0,v}\in f, \tag{20} 
\end{equation}
\begin{equation}          ((\mb ta)_1p)_1= ((\mb t b)_1q)_1\fo \op(v,p^{ia_0},z). \tag{21}  
\end{equation}
Write $\pair{\breve p,\breve q,v}$ for $\pair{((\mb ta)_1p)_0, ((\mb t b)_1q)_0,v}$. Now,  $\mb q_2a=\mb q_2b\fo \vartheta_2(f)$. It then follows from (20) that  
\begin{equation} \mb q_2a\breve p=\mb q_2b\breve q\fo \forall u\, \forall z\, (\op(v,u,z)\imp z\in X_\delta). \tag{22} 
\end{equation}
By (21) and (22), we  get 
$\mb q_2a\breve p((\mb ta)_1p)_1=\mb q_2b\breve q((\mb t b)_1q)_1 \fo z\in X_\delta$, and so we can   construct a closed term $\mb f$ such that 
\[   (\mb fap)_0\sim_\delta (\mb fbq)_0 \ \ \text{and} \ \ (\mb fap)_1= (\mb fbq)_1\fo z=(\mb fap)_0^\delta.\]
The foregoing discussion also shows that, if we let
$(\mb ra)_0=\mb p a_0(\lambda p.(\mb fap)_0)$,
then $c= (\mb ra)_0\sim_\delta (\mb rb)_0=d$. We then wish to find $\mb r$  such that $(\mb ra)_1=(\mb rb)_1\fo y=c^\delta$. For this one needs to  show that $\mb ha=\mb hb\fo f=f_c$ for some closed term $\mb h$. Since this is similar to the proof of the right to left direction of (iii), we leave the details  to the reader.  
\end{proof}


\begin{proof}[Proof of Theorem \ref{canonical}]

We have to  provide closed terms $\mb i$ and $\mb e$ such that
\[  c=d\fo X\in \Prod\Sum\I \ \ \text{implies} \ \ \sigma=\mb ic\sim\mb id\ \ \text{and}\ \  \mb ec=\mb e d\fo X=X_\sigma, \]
for every $X\in\va$. Same task for $\Prod\Sum\W\I$-sets.\\

Let  $\vartheta(u,Y)$ be as in Theorem \ref{def}. Thus $\Prod\Sum\I=\{Y\mid \exists u\, \vartheta(u,Y)\}$. We want $\mb i$ and $\mb e$ such that for every $u,Y\in\va$, 
\[   c=d\fo \vartheta(u,Y)\ \ \text{implies}\ \ \sigma=\mb ic\sim \mb id\ \ \text{and}\ \ \mb ec=\mb e d\fo Y=X_\sigma. \]
The closed terms $\mb i$ and $\mb e$ are obtained by means of the  recursion theorem for pcas. The verification is done by induction on $u$. First, by soundness, there is a closed term $\mb{t}$ such that $c=d\fo \vartheta(u,Y)$ implies
\[ \mb t c=\mb t d\fo \bigvee_{i=0}^4 \vartheta_{i}(u,Y). \]
By using the definition by cases combinator $\mb d$, we can find closed terms $\mb{c}$ and $\mb t_i$ for $i<4$ such that if $c=d\fo  \vartheta(u,Y)$, then for some $i<5$, 

\[ \mb{c} c=\mb{c} d=\bar i\ \ \text{and}\ \ \mb t_ic=\mb t_id\fo \vartheta_{i}(u,Y). \]
\vspace{1mm}	

Case 0. If $\mb c c=\mb 0$, then $\mb t_0 c=\mb t_0 d\fo Y\in\omega$. By Theorem \ref{realizers}, we have a closed term $\mb o$ such that $\mb o\fo X_{\sf N}=\omega$. By soundness, there is a closed term $\mb v_0$ such that $\mb v_0c=\mb v_0d\fo Y\in X_{\sf N}$. This means, that for some $n\in\omega$,
$(\mb v_0c)_0=(\mb v_0d)_0=\bar n$ and 
$(\mb v_0c)_1=(\mb v_0d)_1\fo Y=\dot n$.
Now, ${\sf N}_n=\mb p \mb 0\bar n$ and $X_{{\sf N}_n}=\dot n$. Therefore, in this case, meaning $\mb c c=\mb 0$, we look for $\mb i$ and $\mb e$ such that 
$\mb ic=\mb p\mb 0(\mb v_0c)_0$ and $\mb ec=(\mb v_0c)_1$.
\vspace{1mm}

Case 1. If $\mb cc=\mb 1$, then $\mb t_1c=\mb t_1d\fo Y=\omega$. As before, by soundness, there is a closed term $\mb v_1$ such that $\mb v_1c=\mb v_1d\fo Y=X_{\sf N}$. Therefore, in this case, meaning $\mb cc=\mb 1$, we look for $\mb i$ and $\mb e$ such that 
$\mb ic=o$ and   $\mb ec=\mb v_1c$.
\vspace{1mm}

Case 2. Suppose $\mb c c=\mb 2$. We are in the case where 
\begin{align*}    
	\mb t_2 c=\mb t_2 d&\fo   \exists v\in u\, \vartheta(v,X)\land  \fun(F) \land \dom(F)=X\land \\ \tag{1}
	& \qquad \forall x\in X\, \exists v\in u\, \vartheta(v,F(x))\land Y=\Prod(X,F). 
\end{align*}
for some $X,F\in\va$. By undoing (1), we get that  for  some $\pair{\mathring{c},\mathring{d},v}\in u$, 
$(\mb t_2 c)_{01}=(\mb t_2 d)_{01} \fo \vartheta(v,X)$.
By induction, we can assume 
\begin{equation}  \sigma=\mb i (\mb t_2 c)_{01}\sim\mb i (\mb t_2 d)_{01}=\tau, \tag{2}
\end{equation}
\begin{equation}   \mb e(\mb t_2 c)_{01}=\mb e (\mb t_2 d)_{01} \fo X=X_\sigma.  \tag{3}
\end{equation}
By the properties of equality, on account of (1) and (3), we can find closed terms $\mb q_0,\mb q_1,\mb q_2$ such that
\begin{align*} \tag{4} \mb q_0\mb e c=\mb q_0\mb ed & \fo \fun(F)\land \dom(F)=X_\sigma, \\
	\tag{5} \mb q_1\mb ec=\mb q_1\mb ed &\fo \forall x\in X_\sigma\, \exists v\in u\, \vartheta(v,F(x)), \\
	\tag{6} \mb q_2\mb ec=\mb q_2\mb ed &\fo Y=\Prod(X_\sigma,F). 
\end{align*}

Suppose $a\sim_\sigma b$. By using (4), we can obtain a closed term $\mb r$ such that 
\[  \mb r\mb ec=\mb r\mb ed\fo \forall x\in X_\sigma\, \exists Z\,  \exists \mathring{X}\, (Z\in F\land \op(Z,x,\mathring{X})). \] Therefore, 
there are $Z,\mathring{X}\in\va$ such that 
$\mb r\mb eca=\mb r\mb edb\fo Z\in F\land \op(Z,a^\sigma,\mathring{X})$.
On the other hand, it follows from (5) that
\[ \mb q_1\mb eca=\mb q_1\mb edb\fo \exists v\in u\, \forall Z\, \forall \mathring{X}\, (Z\in F\land \op(Z,a^\sigma,\mathring{X})\imp \vartheta(v,\mathring X)). \]
In particular, for some $\pair{\mathring c,\mathring d,v}\in u$,
$(\mb q_1\mb eca)_1(\mb r\mb eca)=(\mb q_1\mb edb)_1(\mb r\mb edb)\fo  \vartheta(v,\mathring{X})$.
To ease notation, let $r(\mb e,c,a)=(\mb q_1\mb eca)_1(\mb r\mb eca)$. By induction, we can assume that 
$\rho=\mb i r(\mb e,c,a)\sim\mb i r(\mb e,d,b)$ and 
$\mb e r(\mb e,c,a)=\mb e r(\mb e,d,b)\fo \mathring X=X_\rho$.
We have just shown that if $a\sim_\sigma b$ then $ia\sim jb$, where $i=\lambda a.\mb i r(\mb e,c,a)$ and $j=\lambda b.\mb i r(\mb e,d,b)$. By (2),  we thus have $\Prod_\sigma i\sim \Prod_\tau j$.

Let $\mb{prod}$ be as in Theorem \ref{realizers}. That is,  $\mb{prod}\fo X_{\Prod_\sigma i}=\Prod(X_\sigma,F_{\sigma,i})$.
We aim to show that for some closed term $\mb f$, $\mb f \mb ec=\mb f \mb ed\fo F=F_{\sigma,i}$.
If we do so, owing to the soundness, we can cook up, by using (6) and $\mb{prod}$,
a closed term $\mb v_2$ such that
$\mb v_2\mb ec=\mb v_2\mb ed\fo Y=X_{\Prod_\sigma i}$.
In this case, meaning $\mb c c=\mb 2$, we  thus want $\mb i$ and $\mb e$ such that
\[ \mb ic=\Prod_{\mb i (\mb t_2 c)_{01}}(\lambda a.\mb i (\mb q_1\mb eca)_1(\mb r\mb eca))\ \ \text{and}\ \  \mb ec=\mb v_2\mb ec. \]

Let us show how to find $\mb f$. By Theorem \ref{realizers}, $\mb{fun}\fo \fun(F_{\sigma,i})\land \dom(F)=X_\sigma$.  
By (4), it then suffices to look for a closed term $\mb{h}$ such that
\[  \mb h \mb ec=\mb h \mb ed \fo  \forall Z\in F_{\sigma,i}\, (Z\in F). \]
Let $a\sim_\sigma b$. We aim for 
$\mb h \mb eca=\mb h \mb edb\fo \vpair{a^\sigma,X_{\rho}}\in F$,
where $\rho=i a$. As before, we know that for some $Z,\mathring X\in\va$,
$\mb r\mb eca=\mb r\mb edb  \fo Z\in F\land \op(Z,a^\sigma,\mathring{X})$ and
$\mb e r(\mb e,c,a)=\mb e r(\mb e,d,b) \fo \mathring X=X_\rho$.  
By the properties of equality and pairing, we can easily arrange for such  $\mb h$.\vspace{1mm}

Case 3. The case $\mb c c=\mb 3$ is treated similarly.\vspace{1mm}

Case 4. Suppose $\mb cc=\mb 4$. In this case, 
\begin{align*}
	\tag{7}	\mb t_4 c=\mb t_4 d&\fo \exists v\in u\, \vartheta(v,X)\land \exists x\in X\, \exists y\in X\, (Y=\I(x,y)),
\end{align*}
for some $X\in \va$. As in case 2, we can assume 
\[ \sigma=\mb i (\mb t_4 c)_{01}\sim\mb i (\mb t_4 d)_{01}=\tau, \]  
\begin{equation} \mb e(\mb t_4 c)_{01}=\mb e (\mb t_4 d)_{01} \fo X=X_\sigma. \tag{8}  
\end{equation}
By the properties of equality, owing to  (7) and (8), we thus have a closed term $\mb q$ such that 
$\mb q\mb ec=\mb q\mb ed\fo \exists x\in X_\sigma\, \exists y\in X_\sigma\, (Y=\I(x,y))$.
By unfolding the definition, we obtain
\[   a=(\mb q\mb ec)_0\sim_\sigma (\mb qed)_0=\breve a, \]
\[   b=(\mb q\mb ec)_{10}\sim_\sigma (\mb q\mb ed)_{10}=\breve b, \]
\[  (\mb q\mb ec)_{11}=(\mb q\mb ed)_{11}\fo Y=\I(a^\sigma,b^\sigma). \]
By using $\mb{id}\fo X_{\I_\sigma(a,b)}=\I(a^\sigma,b^\sigma)$ from Theorem \ref{realizers}, we then obtain by soundness a closed term $\mb v_4$ such that 
$\mb v_4\mb ec=\mb v_4\mb ed\fo Y=X_{\I_\sigma(a,b)}$. 
In this case, meaning $\mb c c=\mb 4$, we thus want $\mb i$ and $\mb e$ such that
\[ \mb ic=\I_{\mb i (\mb t_4 c)_{01}}((\mb q\mb ec)_0,(\mb q\mb ec)_{10})\ \ \text{and}\ \  \mb ec=\mb v_4\mb ec. \]
\vspace{1mm}

In the $\Prod\Sum\W\I$ case, one considers  $\vartheta^*(u,Y)$  as in Theorem \ref{def}. Thus $\Prod\Sum\W\I=\{Y\mid \exists u\, \vartheta^*(u,Y)\}$. The proof then  proceeds exactly as before, except that now we also need to deal with $\W$-types.  On the other hand, the treatment of $\W$-types is completely analogous to that of $\Prod$-types and $\Sum$-types. The reader can easily fill in the gaps using case 2 as template.
\end{proof}


\begin{thebibliography}{10}
	
	\bibitem{A78}
	P.~Aczel.
	\newblock The type theoretic interpretation of constructive set theory.
	\newblock In {\em Studies in Logic and the Foundations of Mathematics},
	volume~96, pages 55--66. Elsevier, 1978.
	
	\bibitem{A82}
	P.~Aczel.
	\newblock The type theoretic interpretation of constructive set theory: choice
	principles.
	\newblock In {\em Studies in Logic and the Foundations of Mathematics}, volume
	110, pages 1--40. Elsevier, 1982.
	
	\bibitem{A86}
	P.~Aczel.
	\newblock The type theoretic interpretation of constructive set theory:
	inductive definitions.
	\newblock In {\em Logic, methodology and philosophy of science, {VII}
		({S}alzburg, 1983)}, volume 114 of {\em Stud. Logic Found. Math.}, pages
	17--49. North-Holland, Amsterdam, 1986.
	
	\bibitem{czf}
	P.~Aczel and M.~Rathjen.
	\newblock {\em Notes on constructive set theory}.
	\newblock 2010.
	\newblock Available at \url{http://www1.maths.leeds.ac.uk/~rathjen/book.pdf}.
	
	\bibitem{B79}
	M.~J. Beeson.
	\newblock Continuity in intuitionistic set theories.
	\newblock In {\em Studies in Logic and the Foundations of Mathematics},
	volume~97, pages 1--52. Elsevier, 1979.
	
	\bibitem{Beeson79}
	M.~J. Beeson.
	\newblock Goodman's theorem and beyond.
	\newblock {\em Pacific J. Math.}, 84(1):1--16, 1979.
	
	\bibitem{B85}
	M.~J. Beeson.
	\newblock {\em Foundations of constructive mathematics}, volume~6 of {\em
		Ergebnisse der Mathematik und ihrer Grenzgebiete (3) [Results in Mathematics
		and Related Areas (3)]}.
	\newblock Springer-Verlag, Berlin, 1985.
	
	\bibitem{B67}
	E.~Bishop.
	\newblock {\em Foundations of constructive analysis}, volume~60.
	\newblock McGraw-Hill New York, 1967.
	
	\bibitem{Curry30}
	H.~B. Curry.
	\newblock Grundlagen der kombinatorischen {L}ogik.
	\newblock {\em American Journal of Mathematics}, 51:363--384, 1930.
	
	\bibitem{D75}
	R.~Diaconescu.
	\newblock Axiom of choice and complementation.
	\newblock {\em Proceedings of the American Mathematical Society},
	51(1):176--178, 1975.
	
	\bibitem{F75}
	S.~Feferman.
	\newblock A language and axioms for explicit mathematics.
	\newblock In {\em Algebra and logic}, pages 87--139. Springer, 1975.
	
	\bibitem{F79}
	S.~Feferman.
	\newblock Constructive theories of functions and classes.
	\newblock In {\em Logic {C}olloquium '78}, Stud. Logic Found. Math., pages
	159--224. North-Holland, Amsterdam, 1979.
	
	\bibitem{F73}
	H.~Friedman.
	\newblock Some applications of {K}leene's methods for intuitionistic systems.
	\newblock In {\em Cambridge summer school in mathematical logic}, pages
	113--170. Springer, 1973.
	
	\bibitem{friedman_scedrov84}
	H.~Friedman and A.~\v{S}\v{c}edrov.
	\newblock Large sets in intuitionistic set theory.
	\newblock {\em Annals of Pure and Applied Logic}, 27:1--24, 1984.
	
	\bibitem{F19}
	E.~Frittaion.
	\newblock On {G}oodman realizability.
	\newblock {\em Notre Dame J. Formal Logic}, 60(3):523--550, 08 2019.
	
	\bibitem{FR21}
	E.~Frittaion and M.~Rathjen.
	\newblock {Extensional realizability for intuitionistic set theory}.
	\newblock {\em Journal of Logic and Computation}, 31(2):630--653, 3 2021.
	
	\bibitem{G76}
	N.~D. Goodman.
	\newblock The theory of the {G}\"{o}del functionals.
	\newblock {\em The Journal of Symbolic Logic}, 41(3):574--582, 1976.
	
	\bibitem{G78}
	N.~D. Goodman.
	\newblock Relativized realizability in intuitionistic arithmetic of all finite
	types.
	\newblock {\em J. Symbolic Logic}, 43(1):23--44, 1978.
	
	\bibitem{gordeev}
	L.~Gordeev.
	\newblock Proof-theoretical analysis of weak systems of functions and classes.
	\newblock {\em Annals of Pure and Applied Logic}, 38:1--121, 1988.
	
	\bibitem{KTr70}
	G.~Kreisel and A.~S. Troelstra.
	\newblock Formal systems for some branches of intuitionistic analysis.
	\newblock {\em Annals of Mathematical Logic}, 1:229--387, 1970.
	
	\bibitem{M84}
	D.~C. McCarty.
	\newblock {\em Realizability and recursive mathematics}.
	\newblock Thesis (Ph.D.)--The University of Edinburgh, 1985.
	
	\bibitem{M86}
	D.~C. McCarty.
	\newblock Realizability and recursive set theory.
	\newblock {\em Ann. Pure Appl. Logic}, 32(2):153--183, 1986.
	
	\bibitem{M75}
	J.~Myhill.
	\newblock Constructive set theory.
	\newblock {\em The Journal of Symbolic Logic}, 40(3):347--382, 1975.
	
	\bibitem{R05a}
	M.~Rathjen.
	\newblock Constructive set theory and {B}rouwerian principles.
	\newblock {\em J.UCS}, 11(12):2008--2033, 2005.
	
	\bibitem{R06b}
	M.~Rathjen.
	\newblock Choice principles in constructive and classical set theories.
	\newblock In {\em Logic {C}olloquium '02}, volume~27 of {\em Lect. Notes Log.},
	pages 299--326. Assoc. Symbol. Logic, La Jolla, CA, 2006.
	
	\bibitem{R03m}
	M.~Rathjen.
	\newblock The formulae-as-classes interpretation of constructive set theory.
	\newblock In {\em Proof technology and computation}, volume 200 of {\em NATO
		Sci. Ser. III Comput. Syst. Sci.}, pages 279--322. IOS, Amsterdam, 2006.
	
	\bibitem{R06}
	M.~Rathjen.
	\newblock Realizability for constructive {Z}ermelo-{F}raenkel set theory.
	\newblock In {\em Logic {C}olloquium '03}, volume~24 of {\em Lect. Notes Log.},
	pages 282--314. Assoc. Symbol. Logic, La Jolla, CA, 2006.
	
	\bibitem{Schoen24}
	M.~Sch\"onfinkel.
	\newblock {\"U}ber die {B}austeine der mathematischen {L}ogik.
	\newblock {\em Mathematische Annalen}, 92:305--316, 1924.
	
	\bibitem{T73}
	A.~S. Troelstra.
	\newblock {\em Metamathematical investigation of intuitionistic arithmetic and
		analysis}.
	\newblock Lecture Notes in Mathematics, Vol. 344. Springer-Verlag, Berlin,
	1973.
	
	\bibitem{T98}
	A.~S. Troelstra.
	\newblock Realizability.
	\newblock In {\em Handbook of proof theory}, volume 137 of {\em Stud. Logic
		Found. Math.}, pages 407--473. North-Holland, Amsterdam, 1998.
	
	\bibitem{TD88}
	A.~S. Troelstra and D.~van Dalen.
	\newblock {\em Constructivism in mathematics. {V}ol. {II}}, volume 123 of {\em
		Studies in Logic and the Foundations of Mathematics}.
	\newblock North-Holland Publishing Co., Amsterdam, 1988.
	\newblock An introduction.
	
	\bibitem{BS18}
	B.~van~den Berg and L.~van Slooten.
	\newblock Arithmetical conservation results.
	\newblock {\em Indag. Math. (N.S.)}, 29(1):260--275, 2018.
	
	\bibitem{Oosten97}
	J.~van Oosten.
	\newblock Extensional realizability.
	\newblock {\em Annals of Pure and Applied Logic}, 84(3):317 -- 349, 1997.
	
	\bibitem{Oosten08}
	J.~van Oosten.
	\newblock {\em Realizability: an introduction to its categorical side}, volume
	152 of {\em Studies in Logic and the Foundations of Mathematics}.
	\newblock Elsevier B. V., Amsterdam, 2008.
	
\end{thebibliography}
\end{document}